\documentclass{amsproc}
\usepackage{amssymb}
\usepackage{amsfonts}

\newtheorem{theorem}{Theorem}

\newtheorem{corollary}[theorem]{Corollary}

\newtheorem{definition}[theorem]{Definition}
\newtheorem{example}[theorem]{Example}

\newtheorem{lemma}[theorem]{Lemma}

\newtheorem{proposition}[theorem]{Proposition}
\newtheorem{remark}[theorem]{Remark}

\input{tcilatex}

\begin{document}
\title[The transfer]{{The transfer in mod-$p$ group cohomology between 
\begin{equation*}
\Sigma _{p}\int \Sigma _{p^{n-1}}\text{, }\Sigma _{p^{n-1}}\int \Sigma _{p}%
\text{ and }\ \Sigma _{p^{n}}
\end{equation*}%
}}
\author{Nondas E. Kechagias}
\address{Department of Mathematics\\
University of Ioannina\\
Ioannina 45110\\
Greece}
\email{nkechag@uoi.gr}
\urladdr{http://www.math.uoi.gr/\symbol{126}nondas\_k}
\thanks{This paper is in final form and no version of it will be submitted
for publication elsewhere.}

\begin{abstract}
In this work we compute the induced transfer map: 
\begin{equation*}
\bar{\tau}^{\ast }:\func{Im}\left( res^{\ast }:H^{\ast }\left( G\right)
\rightarrow H^{\ast }\left( V\right) \right) \rightarrow \func{Im}\left(
res^{\ast }:H^{\ast }\left( \Sigma _{p^{n}}\right) \rightarrow H^{\ast
}\left( V\right) \right)
\end{equation*}%
in $\func{mod}p$-cohomology. Here $\Sigma _{p^{n}}$\ is the symmetric group
acting on an $n$-dimensional 
%TCIMACRO{\TeXButton{TeX field}{$\mathbb F_p$} }%
%BeginExpansion
$\mathbb F_p$
%EndExpansion
vector space $V$,\ $G=\Sigma _{p^{n},p}$ a $p$-Sylow subgroup, $\Sigma
_{p^{n-1}}\int \Sigma _{p}$,\ \ or $\Sigma _{p}\int \Sigma _{p^{n-1}}$.\
Some answers are given by natural invariants which are related to certain
parabolic subgroups. We also compute a free module basis for certain rings
of invariants over the classical Dickson algebra. This provides a
computation of the image of the appropriate restriction map. Finally, if $%
\xi :\func{Im}\left( res^{\ast }:H^{\ast }\left( G\right) \rightarrow
H^{\ast }\left( V\right) \right) \rightarrow \func{Im}\left( res^{\ast
}:H^{\ast }\left( \Sigma _{p^{n}}\right) \rightarrow H^{\ast }\left(
V\right) \right) $ is the natural epimorphism, then we prove that $\bar{\tau}%
^{\ast }=\xi $ in the ideal generated by the top Dickson algebra generator.\
\ 
\end{abstract}

\maketitle

\section{Introduction-Results}

Let $H$ be a subgroup of a finite group $G$. There are two important maps in
group cohomology going in the opposite direction: the restriction and
transfer. \ The Weyl subgroup acts on the right in group cohomology and the
inclusion $H\hookrightarrow G$\ induces a map\ \ \ \ 
\begin{equation*}
\left( res_{H}^{G}\right) ^{\ast }:H^{\ast }\left( G\right) \rightarrow
H^{\ast }\left( H\right) ^{W_{G}\left( H\right) }
\end{equation*}%
In other words the image of the restriction map is contained in the $%
W_{G}\left( V\right) $-invariants. The role of classical invariant theory in
determining and analyzing cohomology of finite groups is important.\ 

The inclusion \ $H\hookrightarrow G$\ also induces a transfer map 
\begin{equation*}
tr^{\ast }:H^{\ast }\left( H\right) \rightarrow H^{\ast }\left( G\right)
\end{equation*}%
The transfer map plays a fundamental role in group cohomology.

In this work we compute the maps above for particular cases. Some answers
are given by particular invariants which are of the form: a free module
basis over the fundamental object in modular invariant theory, i.e. the
Dickson algebra.

We studied the case $G=%
%TCIMACRO{\TeXButton{TeX field}{\mathbb Z}}%
%BeginExpansion
\mathbb Z%
%EndExpansion
/p\smallint ...\smallint 
%TCIMACRO{\TeXButton{TeX field}{\mathbb Z}}%
%BeginExpansion
\mathbb Z%
%EndExpansion
/p$\ in \cite{Kech1}. We extend those results for $G=\Sigma _{p}\smallint
\Sigma _{p^{n-1}}$\ and $\Sigma _{p^{n-1}}\smallint \Sigma _{p}$. The
methods applied in \cite{Kech1} can not be applied in this case. We compute
the image of the restriction map in Theorem \ref{Restr-Parab} \ for $\Sigma
_{p^{n_{l}}}\smallint ...\smallint \Sigma _{p^{n_{1}}}$. To compute the
transfer, we need to express the previous ring as a module over the Dickson
algebra. We do so in proposition \ref{d'}\ and Theorem \ref{mimic}. Finally,
we show that the induced transfer coincides with the natural, so called,
epimorphism on a certain ideal in Theorems \ref{Main}\ and \ref{final}. \ \ 

Let $V\cong 
%TCIMACRO{\TeXButton{TeX field}{\mathbb F_{p}^{n}}}%
%BeginExpansion
\mathbb F_{p}^{n}%
%EndExpansion
$ be an $n$-dimensional $%
%TCIMACRO{\TeXButton{TeX field}{\mathbb F_p}}%
%BeginExpansion
\mathbb F_p%
%EndExpansion
$\ vector space. Let $\Sigma _{p^{n}}$\ denote the permutations on $V$. Now $%
V$\ has a left action on itself and defines an inclusion: $V\hookrightarrow
\Sigma _{p^{n}}$. Let $\Sigma _{p}\smallint \Sigma _{p^{n-1}}$\ denote the
semidirect product of $\Sigma _{p}$ with $\left( \Sigma _{p^{n-1}}\right)
^{p}$\ with $\Sigma _{p}$\ acting by permuting factors. And for $\Sigma
_{p^{n-1}}\smallint \Sigma _{p}$\ respectively. Let 
\begin{equation*}
\Sigma _{p^{n},p}:=%
%TCIMACRO{\TeXButton{TeX field}{\mathbb Z}}%
%BeginExpansion
\mathbb Z%
%EndExpansion
/p\smallint ...\smallint 
%TCIMACRO{\TeXButton{TeX field}{\mathbb Z}}%
%BeginExpansion
\mathbb Z%
%EndExpansion
/p
\end{equation*}%
which is a $p$-Sylow subgroup of $\Sigma _{p^{n},p}$ containing $V$. The
maximal elementary abelian $p$-subgroup $V$\ is contained by both $\Sigma
_{p}\smallint \Sigma _{p^{n-1}}$\ and $\Sigma _{p^{n-1}}\smallint \Sigma
_{p} $. \ 

Simple coefficients are taken in $%
%TCIMACRO{\TeXButton{TeX field}{\mathbb F_p}}%
%BeginExpansion
\mathbb F_p%
%EndExpansion
\cong 
%TCIMACRO{\TeXButton{TeX field}{\mathbb Z}}%
%BeginExpansion
\mathbb Z%
%EndExpansion
/p$ where $p$ is an odd prime. For $p=2$ minor modifications are needed and
left to the interested reader. Hence $H^{\ast }\left( G\right) $\ stands for 
$H^{\ast }\left( G,%
%TCIMACRO{\TeXButton{TeX field}{\mathbb Z}}%
%BeginExpansion
\mathbb Z%
%EndExpansion
/p\right) $. \ 

It is known that 
\begin{equation*}
H^{\ast }\left( V\right) \cong \QATOPD\{ . {%
%TCIMACRO{\TeXButton{TeX field}{\mathbb F_p}}%
%BeginExpansion
\mathbb F_p%
%EndExpansion
[y_{1},\cdots ,y_{n}]\text{, for }p=2}{E_{%
%TCIMACRO{\TeXButton{TeX field}{\mathbb F_p}}%
%BeginExpansion
\mathbb F_p%
%EndExpansion
}\left( x_{1},\cdots x_{n}\right) \otimes 
%TCIMACRO{\TeXButton{TeX field}{\mathbb F_p}}%
%BeginExpansion
\mathbb F_p%
%EndExpansion
[y_{1},\cdots ,y_{n}]}
\end{equation*}%
It is known that the Weyl subgroups $W_{\Sigma _{p^{n}}}\left( V\right) $, $%
W_{\Sigma _{p^{n},p}}\left( V\right) $, $W_{\Sigma _{p}\smallint \Sigma
_{p^{n-1}}}\left( V\right) $\ and $W_{\Sigma _{p^{n-1}}\smallint \Sigma
_{p}}\left( V\right) $\ are the general linear group $GL(n,\mathbb{F}_{p})$,
the upper triangular subgroup $U_{n}$, and the parabolic subgroups $P\left(
1,n-1\right) $\ and $P\left( n-1,1\right) $\ respectively. Here {\small {%
\begin{equation*}
P\left( k,n-k\right) =\left\{ \left( 
\begin{array}{cc}
A & C \\ 
\mathbf{0} & B%
\end{array}%
\right) |\;A\in GL(k),B\in GL(n-k))\right\}
\end{equation*}%
}}

Kuhn (\cite{Kuhn}) proved that the following diagram is commutative and this
is the key point for our study:%
\begin{equation*}
\begin{array}{ccc}
H^{\ast }\left( \Sigma _{p}\smallint \Sigma _{p^{n-1}}\right) & \overset{%
tr^{\ast }}{\rightarrow } & H^{\ast }\left( \Sigma _{p^{n}}\right) \\ 
\downarrow \left( res_{V}^{\Sigma _{p}\smallint \Sigma _{p^{n-1}}}\right)
^{\ast } &  & \downarrow \left( res_{V}^{\Sigma _{p^{n}}}\right) ^{\ast } \\ 
H^{\ast }\left( V\right) ^{W_{\Sigma _{p}\smallint \Sigma _{p^{n-1}}}\left(
V\right) } & \overset{\tau ^{\ast }}{\rightarrow } & H^{\ast }\left(
V\right) ^{W_{\Sigma _{p^{n}}}\left( V\right) }%
\end{array}%
\end{equation*}

In this work we investigate the induced transfer homomorphisms: 
\begin{equation*}
\func{Im}\left( res^{\ast }:H^{\ast }\left( G\right) \rightarrow H^{\ast
}\left( V\right) \right) \overset{\bar{\tau}^{\ast }}{\rightarrow }\func{Im}%
\left( res^{\ast }:H^{\ast }\left( \Sigma _{p^{n}}\right) \rightarrow
H^{\ast }\left( V\right) \right)
\end{equation*}%
For $G=\Sigma _{p^{n},p}$, $\Sigma _{p}\int \Sigma _{p^{n-1}}$ and $\Sigma
_{p^{n-1}}\smallint \Sigma _{p}$. The problem reduces to find free module
bases for certain algebras\ of modular invariants. This is a hard problem
for a general parabolic subgroup.

The restriction map is not an onto map and our first task is to compute its
image. Please note that for $p=2$\ the restriction map is onto. We give an
invariant theoretic proof of the following Theorem first proved by Mui (\cite%
{Mui}) using cohomological methods in section 3. It requires technical
results from group cohomology and invariant theory.

\textbf{Theorem \cite{Mui}} \textbf{\ref{Steen-May}} \textit{The image }$%
\func{Im}\left( res^{\ast }:H^{\ast }\left( \Sigma _{p^{n},p}\right)
\rightarrow H^{\ast }\left( V\right) \right) $\textit{\ is isomorphic with
the tensor product between an exterior and a polynomial algebra }%
\begin{equation*}
E_{%
%TCIMACRO{\TeXButton{TeX field}{\mathbb F_p}}%
%BeginExpansion
\mathbb F_p%
%EndExpansion
}\left( \hat{M}_{1,0},\hat{M}_{2,1}\hat{L}_{1}^{\left( p-3\right) /2},\cdots
,\hat{M}_{n,n-1}\hat{L}_{n-1}^{\left( p-3\right) /2}\right) \otimes H_{n}^{t}
\end{equation*}%
Definitions and notation are given in section 2.

Since the transfer is an additive map (and the identity on the Dickson
algebra), it is important to describe these images of the appropriate rings
as modules over the Dickson algebra ($H^{\ast }\left( V\right) ^{GL(n,%
\mathbb{F}_{p})}$). The bulk of this work is to that direction.

As an application of last Theorem we derive the next proposition in section
3. The image is given by natural invariants which have the following form.

\textbf{Proposition} \textbf{\ref{decomp-restric}} \textit{The image }$\func{%
Im}\left( res^{\ast }:H^{\ast }\left( \Sigma _{p^{n},p}\right) \rightarrow
H^{\ast }\left( V\right) \right) $\textit{\ is isomorphic with }%
\begin{equation*}
H_{n}^{t}\tbigoplus\limits_{i}\tbigoplus\limits_{s_{i}}H_{n}^{t}\hat{M}%
_{i,s_{1},...,s_{k-1},i-1}\hat{L}_{i-1}^{\left( p-3\right)
/2}\tprod\limits_{1}^{k-1}\hat{L}_{s_{t}}^{\left( p-3\right)
/2}\tprod\limits_{1}^{k-1}\hat{L}_{s_{t}+1}^{\left( p-3\right) /2}
\end{equation*}%
\textit{Here }$k\leq i\leq n$\textit{\ and }$0\leq s_{1}<...<s_{k-1}<i-1$%
\textit{.}

Let $I=(n_{l},...,n_{1})$ be a sequence of positive integers such that $%
\tsum n_{i}=n$ and $P\left( I\right) $\ the associated parabolic subgroup.
We call 
\begin{equation*}
D_{n}:=\left( 
%TCIMACRO{\TeXButton{TeX field}{\mathbb F_p}}%
%BeginExpansion
\mathbb F_p%
%EndExpansion
[y_{1},...,y_{n}]\right) ^{GL(n,\mathbb{F}_{p})}
\end{equation*}%
(the classical Dickson algebra) and 
\begin{equation*}
%TCIMACRO{\TeXButton{TeX field}{\mathbb F_p}}%
%BeginExpansion
\mathbb F_p%
%EndExpansion
\left( I\right) :=\left( 
%TCIMACRO{\TeXButton{TeX field}{\mathbb F_p}}%
%BeginExpansion
\mathbb F_p%
%EndExpansion
[y_{1},\cdots ,y_{n}]\right) ^{P\left( I\right) }
\end{equation*}%
Implementing last Theorem and the ring $H^{\ast }\left( V\right) ^{P\left(
I\right) }$, we compute the image of the restriction map in section 3.

\textbf{Theorem \ref{Restr-Parab} }\textit{The image\ }$\func{Im}\left(
res^{\ast }:H^{\ast }\left( \Sigma _{p^{n_{l}}}\smallint ...\smallint \Sigma
_{p^{n_{1}}}\right) \rightarrow H^{\ast }\left( V\right) \right) $\textit{\
is isomorphic to the subalgebra generated by }%
\begin{equation*}
\left\{ 
\begin{array}{c}
\hat{d}_{\nu _{i},\nu _{i}-k_{i}}\ ,\hat{M}_{\nu _{i},\nu _{i}-k_{i}}\left( 
\hat{L}_{\nu _{i}}\right) ^{p-2},\hat{M}_{\nu _{i},\nu _{i}-k_{j},\nu
_{i}-k_{i}}\left( \hat{L}_{\nu _{i}}\right) ^{p-2}|\  \\ 
1\leq i\leq \ell ,\ 1\leq k_{i}\leq n_{i},k_{i}<k_{j}<\nu _{i},\ \nu
_{i}=\sum\limits_{t=1}^{i}n_{t}%
\end{array}%
\right\}
\end{equation*}%
\textit{along with certain relations. }

For notation and relations between the generators please see Theorem \ref%
{invar-Bock-Parab} in section 2.

It is a hard problem to express the subalgebra above as a free module over
the appropriate subalgebra of the Dickson algebra. Instead we study certain
rings of invariants of parabolic subgroups.

It is known that $%
%TCIMACRO{\TeXButton{TeX field}{\mathbb F_p}}%
%BeginExpansion
\mathbb F_p%
%EndExpansion
\left( I\right) $\ is a finitely generated free module over $D_{n}$. In
order to provide a free basis, we define a new generating set for $%
%TCIMACRO{\TeXButton{TeX field}{\mathbb F_p}}%
%BeginExpansion
\mathbb F_p%
%EndExpansion
\left( 1,n-1\right) $ and $%
%TCIMACRO{\TeXButton{TeX field}{\mathbb F_p}}%
%BeginExpansion
\mathbb F_p%
%EndExpansion
\left( n-1,1\right) $. There are two advantages for this new set. Mainly, it
is closed under the action of Steenrod's algebra and secondly the algebra
generators for $D_{n}$ can be decomposed with respect to the new ones. We
prove the following proposition in section 4.\ 

\textbf{Proposition} \textbf{\ref{d'}} \textit{Let }$I=\left( 1,n-1\right) $%
\textit{, then }%
\begin{equation*}
%TCIMACRO{\TeXButton{TeX field}{\mathbb F_p}}%
%BeginExpansion
\mathbb F_p%
%EndExpansion
\left( I\right) =%
%TCIMACRO{\TeXButton{TeX field}{\mathbb F_p}}%
%BeginExpansion
\mathbb F_p%
%EndExpansion
[h_{1}^{p-1},d_{n,i}\left( I\right) \ |\ 1\leq i\leq n-1]
\end{equation*}%
\begin{equation*}
H^{\ast }\left( V\right) ^{P\left( I\right) }\cong 
%TCIMACRO{\TeXButton{TeX field}{\mathbb F_p}}%
%BeginExpansion
\mathbb F_p%
%EndExpansion
\left( I\right) \oplus 
%TCIMACRO{\TeXButton{TeX field}{\mathbb F_p}}%
%BeginExpansion
\mathbb F_p%
%EndExpansion
\left( I\right) \left[ M_{1,0}h_{1}^{p-2}\tbigoplus%
\limits_{t_{i}}M_{n,t_{1},...,t_{k}}L_{n}^{p-2}\right]
\end{equation*}%
\textit{Here }$1\leq t_{k}$\textit{\ and\ }$0\leq t_{1}<...<t_{k}\leq n-1$%
\textit{.}

Kuhn and Mitchell described $%
%TCIMACRO{\TeXButton{TeX field}{\mathbb F_p}}%
%BeginExpansion
\mathbb F_p%
%EndExpansion
\left( I\right) $\ using appropriate Dickson algebra generators in \cite%
{Kuhn-Mitc}.\ Their set is elegant and more easily described than ours, but
their set is not closed under the action of Steenrod's algebra, and their
set is not as useful as ours is in computations.

The next Theorem provides a free module basis for $%
%TCIMACRO{\TeXButton{TeX field}{\mathbb F_p}}%
%BeginExpansion
\mathbb F_p%
%EndExpansion
\left( n-1,1\right) $ over $D_{n}$ proved in section 5.\ 

For each $t$, $1\leq t\leq n-1$, we define the set of all $\left( n-t\right) 
$-tuples\ \ 
\begin{equation*}
\mathcal{M}(n-2,t)=\{M=(p,m_{t},...,m_{n-2})\ |\ 0\leq m_{i}\leq p-1\}
\end{equation*}%
and, for each $M\in \mathcal{M}(n-2,t)$\ we define 
\begin{equation*}
d_{n-1}^{M}=d_{n-1,t-1}^{p}d_{n-1,t}^{m_{t}}...d_{n-1,n-2}^{m_{n-2}}
\end{equation*}

\textbf{Theorem} \textbf{\ref{mimic}} We have 
\begin{equation*}
B_{D_{n}}\left( 
%TCIMACRO{\TeXButton{TeX field}{\mathbb F_p}}%
%BeginExpansion
\mathbb F_p%
%EndExpansion
(n-1,1)\right) =\tbigcup\limits_{t=1}^{n-1}\{d_{n-1}^{M}\,|\,M\in \mathcal{M}%
(n-2,t)\}
\end{equation*}%
\textit{as a free module basis for 
%TCIMACRO{\TeXButton{TeX field}{$\mathbb F_p$}}%
%BeginExpansion
$\mathbb F_p$%
%EndExpansion
}$(n-1,1)$\textit{\ over }$D_{n}$\textit{.\ \ }

The following corollary is the main result in this work.

\textbf{Corollary \ref{mainCorollar} }\textit{i) }$\func{Im}\left( res^{\ast
}:H^{\ast }\left( \Sigma _{p}\int \Sigma _{p^{n-1}}\right) \rightarrow
H^{\ast }\left( V\right) \right) $\textit{\ is isomorphic to a free module
over }$D_{n}$\textit{\ on \ }%
\begin{equation*}
\left\{ 
\begin{array}{c}
\hat{M}_{1,0}\hat{L}_{1}^{\left( p-2\right) }\hat{h}_{1}^{(p-1)m},\hat{M}%
_{n,s_{1},...,s_{k}}\hat{L}_{n}^{\left( p-2\right) }d_{n,0}^{\left( \left[ 
\frac{k+1}{2}\right] -1\right) }\hat{h}_{1}^{(p-1)m}\ |\  \\ 
0\leq m<A_{1},k\leq n,1\leq s_{k},0\leq s_{1}<...<s_{k}\leq n-1%
\end{array}%
\right\}
\end{equation*}%
\textit{Here }$A_{1}=p^{n-1}+...+p$\textit{.}

\textit{ii) }$\func{Im}\left( res^{\ast }:H^{\ast }\left( \Sigma
_{p^{n-1}}\int \Sigma _{p}\right) \rightarrow H^{\ast }\left( V\right)
\right) $\textit{\ is isomorphic to a free module over }$D_{n}$\textit{\ on
\ }%
\begin{equation*}
\left\{ 
\begin{array}{c}
\hat{M}_{n,n-1}\hat{L}_{n}^{\left( p-2\right) }f,\hat{M}%
_{n-1,s_{1},...,s_{k}}\hat{L}_{n-1}^{\left( p-2\right) }d_{n-1,0}^{\left( %
\left[ \frac{k+1}{2}\right] -1\right) }g\ |\  \\ 
f,g\in B_{D_{n}}\left( 
%TCIMACRO{\TeXButton{TeX field}{\mathbb F_p}}%
%BeginExpansion
\mathbb F_p%
%EndExpansion
(n-1,1)\right) ,\text{ }k\leq n-1,0\leq s_{1}<...<s_{k}\leq n-1%
\end{array}%
\right\}
\end{equation*}

Finally, the transfer map is studied in the last section. There is a natural
description of $%
%TCIMACRO{\TeXButton{TeX field}{\mathbb F_p}}%
%BeginExpansion
\mathbb F_p%
%EndExpansion
\left( 1,n-1\right) $ or $%
%TCIMACRO{\TeXButton{TeX field}{\mathbb F_p}}%
%BeginExpansion
\mathbb F_p%
%EndExpansion
\left( n-1,1\right) $ as a polynomial algebra (proposition \textbf{\ref{d'}}
or as described in \cite{Kuhn-Mitc}). According to last corollary, there is
an alternate description of it as a free module over the Dickson algebra.
The natural epimorphisms 
\begin{equation*}
\xi :%
%TCIMACRO{\TeXButton{TeX field}{\mathbb F_p}}%
%BeginExpansion
\mathbb F_p%
%EndExpansion
\left( 1,n-1\right) \rightarrow D_{n}\text{ and }\xi :%
%TCIMACRO{\TeXButton{TeX field}{\mathbb F_p}}%
%BeginExpansion
\mathbb F_p%
%EndExpansion
\left( n-1,1\right) \rightarrow D_{n}
\end{equation*}%
which "rewrites" an element of the polynomial algebra in terms of the free
module basis are shown to be equal with the induced transfer maps. Let us
consider an example.

\textbf{Example} Let $n=3$\ and $p=2$. $%
%TCIMACRO{\TeXButton{TeX field}{\mathbb F_p}}%
%BeginExpansion
\mathbb F_p%
%EndExpansion
\left( 2,1\right) =%
%TCIMACRO{\TeXButton{TeX field}{\mathbb F_p}}%
%BeginExpansion
\mathbb F_p%
%EndExpansion
\left[ d_{2,0},d_{2,1},d_{3,2}\right] $ and the basis is $B=\left\{
d_{2,0}^{i}d_{2,1}^{j},d_{2,0}^{2}d_{2,1}^{j},d_{2,1}^{2}|0\leq i,j\leq
1\right\} $. We need to describe the way in which the three generators of $%
%TCIMACRO{\TeXButton{TeX field}{\mathbb F_p}}%
%BeginExpansion
\mathbb F_p%
%EndExpansion
\left( 2,1\right) $ can be written in terms of $B$\ and $D_{3}$. Here is the
way:\ 
\begin{eqnarray*}
d_{2,0}d_{2,1}^{2} &=&d_{3,0}+d_{3,2}d_{2,0} \\
d_{2,1}^{3} &=&d_{3,1}+d_{3,2}d_{2,1}+d_{2,0}^{2} \\
d_{2,0}^{3} &=&d_{3,1}d_{2,0}+d_{3,0}d_{2,1}
\end{eqnarray*}%
Suppose we want to find $\xi \left( d_{2,0}^{2}d_{2,1}^{7}\right) $.
According to $B$\ and the relations above, this element "rewrites" as follows%
\begin{equation*}
d_{2,0}^{2}d_{2,1}^{7}=d_{3,0}^{2}d_{3,1}+d_{3,0}^{2}d_{3,2}d_{2,1}+d_{3,0}^{2}d_{2,0}^{2}+d_{3,2}^{3}d_{2,0}^{2}d_{2,1}+d_{3,0}d_{3,2}^{2}d_{2,0}d_{2,1}
\end{equation*}%
Thus $\xi \left( d_{2,0}^{2}d_{2,1}^{7}\right) =d_{3,0}^{2}d_{3,1}$.\ 

\textbf{Theorem} \textbf{\ref{Main}} \textit{Let }$\xi :$\textit{\ 
%TCIMACRO{\TeXButton{TeX field}{$\mathbb F_p$}}%
%BeginExpansion
$\mathbb F_p$%
%EndExpansion
}$(n-1,1)\longrightarrow D_{n}$\textit{\ be the natural epimorphism with
respect to the given free module basis }$B$\textit{\ and }$\bar{\tau}^{\ast
}:%
%TCIMACRO{\TeXButton{TeX field}{\mathbb F_p}}%
%BeginExpansion
\mathbb F_p%
%EndExpansion
(n-1,1)\rightarrow D_{n}$\textit{\ the transfer map. Then }$\xi =\bar{\tau}%
^{\ast }$\textit{.}

The advantage of the map $\xi $ is that it calculates $\bar{\tau}^{\ast }$.

Although the transfer map satisfies the nice property described in last
Theorem for the polynomial part of the ring of invariants, it does not for
the exterior part. Please see example \ref{Paradeigma}. But the transfer
coincides with the map $\xi $ in the ideal generated by the top Dickson
algebra generator.

\textbf{Theorem}\ \textbf{\ref{final}} \textit{Let }$\xi ,\bar{\tau}^{\ast }:%
\func{Im}\left( res_{V}^{\Sigma _{p^{n},p}}\right) ^{\ast }\rightarrow \func{%
Im}\left( res_{V}^{\Sigma _{p^{n}}}\right) ^{\ast }$\textit{\ the rewriting
and the induced transfer maps. Then }$\xi =\bar{\tau}^{\ast }$\textit{\ in
the ideal generated by }$\left( d_{n,0}\right) $\textit{. }\ \ 

Our method strongly depends on the action of Steenrod's algebra on the rings
of invariants. This action is the key ingredient in the proof of Theorem \ref%
{Steen-May}\ which is the building block for the computation of the images
of the appropriate restriction maps.\ This method was inspired by a similar
method used by Adem and Milgram VI 1 in \cite{Ad-Mil}. All background
material can be found in this excellent account. For the computation of the
free module bases, we follow Campbell and Hughes \cite{C-H}. Taking into
account proposition \ref{decomp-restric} which is a long and technical
result, the familiar reader may proceed to sections 5 and 6.

We thank the referee and N. Kuhn very much for their suggestions regarding
the exposition of this work.

\section{The rings of invariants}

Let us repeat some classical results from the literature. Let $G=GL(n,%
%TCIMACRO{\TeXButton{TeX field}{\mathbb F_p}}%
%BeginExpansion
\mathbb F_p%
%EndExpansion
)$, $B_{n},$ or $U_{n}$ be the general linear group, the Borel subgroup, and
the upper triangular subgroup with 1's on the diagonal, respectively. $G$
acts as usual on $V$. Let $I=(n_{l},...,n_{1})$ be an ordered sequence of
positive integers such that $\tsum n_{i}=n$. We order such sequences as
above by refinements: $I\leq I^{\prime }$ if $I$ is a refinement of $%
I^{\prime }$. For example $(1,...,1)\leq (n_{1},n_{2})\leq (n)$. Given such
a sequence $I$ let $V^{1}\subset ...\subset V^{l}=V$ be defined by 
\begin{equation*}
V^{i}=<e_{1},e_{2},...,e_{(n_{1}+...+n_{i})}>
\end{equation*}%
This is called a flag by Kuhn \cite{Kuhn}. It is well known that the set 
\begin{equation*}
P(I):=\{g\in GL(n,%
%TCIMACRO{\TeXButton{TeX field}{\mathbb F_p}}%
%BeginExpansion
\mathbb F_p%
%EndExpansion
)\;|\;\forall i\;g(V^{i})=V^{i}\}
\end{equation*}%
\begin{equation*}
P\left( I\right) =\left\{ \left( 
\begin{array}{ccc}
GL_{n_{1}} & \ast & \ast \\ 
0 & \ddots & \ast \\ 
0 & 0 & GL_{n_{\ell }}%
\end{array}%
\right) \right\} \leq GL(n,\mathbb{F}_{p})
\end{equation*}%
is a subgroup of $GL(n,%
%TCIMACRO{\TeXButton{TeX field}{\mathbb F_p}}%
%BeginExpansion
\mathbb F_p%
%EndExpansion
)$ called a parabolic subgroup related to the partition $I$. Moreover, if $G$
is a subgroup of $GL(n,%
%TCIMACRO{\TeXButton{TeX field}{\mathbb F_p}}%
%BeginExpansion
\mathbb F_p%
%EndExpansion
)$ containing the Borel subgroup $B_{n}$, then $G=P(I)$ for some sequence $I$%
, (\cite{Cart} page 112).

Since $H^{\ast }\left( V\right) =E_{%
%TCIMACRO{\TeXButton{TeX field}{\mathbb F_p}}%
%BeginExpansion
\mathbb F_p%
%EndExpansion
}\left( x_{1},\cdots x_{n}\right) \otimes 
%TCIMACRO{\TeXButton{TeX field}{\mathbb F_p}}%
%BeginExpansion
\mathbb F_p%
%EndExpansion
[y_{1},\cdots ,y_{n}]$, the object of study is

\begin{equation*}
\left( E_{%
%TCIMACRO{\TeXButton{TeX field}{\mathbb F_p}}%
%BeginExpansion
\mathbb F_p%
%EndExpansion
}\left( x_{1},\cdots x_{n}\right) \otimes 
%TCIMACRO{\TeXButton{TeX field}{\mathbb F_p}}%
%BeginExpansion
\mathbb F_p%
%EndExpansion
[y_{1},\cdots ,y_{n}]\right) ^{P(I)}
\end{equation*}

The classical Dickson algebra, $D_{n}=\left( 
%TCIMACRO{\TeXButton{TeX field}{\mathbb F_p}}%
%BeginExpansion
\mathbb F_p%
%EndExpansion
[y_{1},\cdots ,y_{n}]\right) ^{GL(n,\mathbb{F}_{p})}$, is described as
follows. Let 
\begin{equation*}
h_{i}=\prod\limits_{v\in V^{i-1}}(y_{i}-v)\text{ and }L_{n}=%
\prod_{1}^{n}h_{i}
\end{equation*}%
Let $L_{n,i}$ be the determinant of the $n\times n$ matrix $\left( 
\begin{array}{ccc}
y_{1} & \cdots & y_{n} \\ 
\vdots & \cdots & \vdots \\ 
y_{1}^{p^{n}} & \cdots & y_{n}^{p^{n}}%
\end{array}%
\right) $ where the $i+1$-row is missing, i.e. the row $\left( 
\begin{array}{ccc}
y_{1}^{p^{i}}, & \cdots , & y_{n}^{p^{i}}%
\end{array}%
\right) $. Moreover, $L_{n}=L_{n,n}$ and $L_{n,0}=L_{n}^{p}$.

Let $L_{n,i}\left( \hat{t}\right) =\det \left( 
\begin{array}{ccccc}
y_{1} & \cdots & \hat{y}_{t} & \cdots & y_{n} \\ 
\vdots & \cdots & \vdots & \cdots & \vdots \\ 
y_{1}^{p^{n-1}} & \cdots & \hat{y}_{t}^{p^{n-1}} & \cdots & y_{n}^{p^{n-1}}%
\end{array}%
\right) $\ where the $i+1$-row is missing.\ Now the following formula holds:%
\begin{equation}
L_{n}=\left( -1\right) ^{t-1}[y_{t}L_{n,n-1}\left( \hat{t}\right)
-y_{t}^{p}L_{n,1}\left( \hat{t}\right) +...+\left( -1\right)
^{n-1}y_{t}^{p^{n-1}}L_{n,n-1}\left( \hat{t}\right) ]  \label{L(1)}
\end{equation}

Finally, let 
\begin{equation*}
d_{n,i}=\frac{L_{n,i}}{L_{n}}
\end{equation*}%
The degrees of the previous elements are $|h_{i}|=2p^{i-1}$, $|L_{n}|=2\frac{%
p^{n}-1}{p-1}$, and $|d_{n,i}|=2\left( p^{n}-p^{i}\right) $.

We shall also need the matrix $\omega $ which consists of $1$'s along the
antidiagonal for the transpose of these groups, please see remark \ref%
{transpose}.\ \ 

\begin{definition}
Let $f\in H^{\ast }\left( V\right) $, then $\hat{f}$\ stands for $\omega f$.
In particular $\hat{h}_{i}=\omega h_{i}$ or $\hat{h}_{i}=\prod\limits_{v\in
\left\langle y_{n+2-i},...,y_{n}\right\rangle }(y_{n+1-i}-v)$. \ \ \ 
\end{definition}

\begin{theorem}[Dickson]
\cite{Dic}$D_{n}=$%
%TCIMACRO{\TeXButton{TeX field}{$\mathbb F_p$}}%
%BeginExpansion
$\mathbb F_p$%
%EndExpansion
$[d_{n,0},\cdots ,d_{n,n-1}]$.
\end{theorem}

\begin{theorem}[Mui]
\cite{Mui}i) $H_{n}:=\left( 
%TCIMACRO{\TeXButton{TeX field}{\mathbb F_p}}%
%BeginExpansion
\mathbb F_p%
%EndExpansion
[y_{1},\cdots ,y_{n}]\right) ^{U_{n}}=$%
%TCIMACRO{\TeXButton{TeX field}{$\mathbb F_p$}}%
%BeginExpansion
$\mathbb F_p$%
%EndExpansion
$[h_{n},\cdots ,h_{1}]$ and 
\begin{equation*}
H_{n}^{t}:=\left( 
%TCIMACRO{\TeXButton{TeX field}{\mathbb F_p}}%
%BeginExpansion
\mathbb F_p%
%EndExpansion
[y_{1},\cdots ,y_{n}]\right) ^{U_{n}^{t}}=%
%TCIMACRO{\TeXButton{TeX field}{\mathbb F_p}}%
%BeginExpansion
\mathbb F_p%
%EndExpansion
[\hat{h}_{n},\cdots ,\hat{h}_{1}]
\end{equation*}

ii) $\left( 
%TCIMACRO{\TeXButton{TeX field}{\mathbb F_p}}%
%BeginExpansion
\mathbb F_p%
%EndExpansion
[y_{1},\cdots ,y_{n}]\right) ^{B_{n}}=$%
%TCIMACRO{\TeXButton{TeX field}{$\mathbb F_p$}}%
%BeginExpansion
$\mathbb F_p$%
%EndExpansion
$[(h_{n})^{p-1},\cdots ,(h_{1})^{p-1}]$ and 
\begin{equation*}
\left( 
%TCIMACRO{\TeXButton{TeX field}{\mathbb F_p}}%
%BeginExpansion
\mathbb F_p%
%EndExpansion
[y_{1},\cdots ,y_{n}]\right) ^{B_{n}^{t}}=%
%TCIMACRO{\TeXButton{TeX field}{\mathbb F_p}}%
%BeginExpansion
\mathbb F_p%
%EndExpansion
[(\hat{h}_{n})^{p-1},\cdots ,(\hat{h}_{1})^{p-1}]
\end{equation*}
\end{theorem}

Relations between the generators of rings of invariants are given as follows:

\begin{proposition}
\cite{Kech1}\label{formula} $d_{n,n-i}=\sum\limits_{1\leq j_{1}<\cdots
<j_{i}\leq n}\prod\limits_{s=1}^{i}\left( h_{j_{s}}^{p-1}\right)
^{p^{n-i+s-j_{s}}}$.
\end{proposition}

\begin{corollary}
\label{corol1}$d_{n,n-i}=d_{n-1,n-i}h_{n}^{p-1}+d_{n-1,n-i-1}^{p}$.
\end{corollary}

\begin{theorem}[Kuhn and Mitchell]
\cite{Kuhn-Mitc} Let $I=(n_{l},\cdots ,n_{1})$. \newline
i) 
%TCIMACRO{\TeXButton{TeX field}{$\mathbb F_p$}}%
%BeginExpansion
$\mathbb F_p$%
%EndExpansion
$(I):=$ 
%TCIMACRO{\TeXButton{TeX field}{$\mathbb F_p$}}%
%BeginExpansion
$\mathbb F_p$%
%EndExpansion
$[d_{\nu _{i},\nu _{i}-k_{i}}\ |\ 1\leq i\leq \ell ,\ 1\leq k_{i}\leq
n_{i},\ \nu _{i}=\sum\limits_{t=1}^{i}n_{t}]$.

ii) 
%TCIMACRO{\TeXButton{TeX field}{$\mathbb F_p$}}%
%BeginExpansion
$\mathbb F_p$%
%EndExpansion
$(I)^{t}:=$ 
%TCIMACRO{\TeXButton{TeX field}{$\mathbb F_p$}}%
%BeginExpansion
$\mathbb F_p$%
%EndExpansion
$[\hat{d}_{\nu _{i},\nu _{i}-k_{i}}\ |\ 1\leq i\leq \ell ,\ 1\leq k_{i}\leq
n_{i},\ \nu _{i}=\sum\limits_{t=1}^{i}n_{t}]$.
\end{theorem}

All the rings of invariants considered in this work are algebras over the
Steenrod algebra. The action of Steenrod's algebra on Dickson algebra
elements has been completely computed in \cite{Kech2}. We repeat here the
following Theorem applied several times in this work.

\begin{theorem}
\label{Steenr-action-dn}(\cite{Kech2} 36, page 170) i) Let $q=\Sigma
_{1}^{n-1}a_{t}p^{t+l}$ such that $p-1\geq a_{t}\geq a_{t-1}>a_{i-1}=0$. Then%
\begin{equation*}
P^{q}d_{n,0}^{p^{l}}=d_{n,0}^{p^{l}}\left( -1\right) ^{a_{n-1}}\Pi _{i}^{n-1}%
\binom{a_{t}}{a_{t-1}}d_{n,t}^{p^{l}\left( a_{t}-a_{t-1}\right) }
\end{equation*}%
Otherwise, $P^{q}d_{n,0}^{p^{l}}=0$.

ii) Let $q=\Sigma _{1}^{n-1}a_{t}p^{t+l}$ such that $p-1\geq a_{t}\geq
a_{t-1}>a_{i}=0$ and $a_{i}+1$\ $\geq a_{i-1}\geq a_{t}\geq a_{t-1}\geq 0$.
Then%
\begin{equation*}
P^{q}d_{n,i}^{p^{l}}=
\end{equation*}%
\begin{equation*}
d_{n,i}^{p^{l}}\left( -1\right) ^{a_{n-1}}\left( \Pi _{i+1}^{n-1}\binom{a_{t}%
}{a_{t-1}}\right) \binom{a_{i}+1}{a_{i-1}}\left( \Pi _{s}^{i-1}\binom{a_{t}}{%
a_{t-1}}\right) \Pi _{s}^{n-1}d_{n,t}^{p^{l}\left( a_{t}-a_{t-1}\right) }
\end{equation*}%
Here $a_{s-1}=0$. Otherwise, $P^{q}d_{n,0}^{p^{l}}=0$.\ \ 
\end{theorem}

We need some technical results for the proof of Theorem \ref{Steen-May}. Let 
\begin{equation}
h_{i}\left( \hat{\jmath}\right) :=\prod\limits_{v\in \left\langle y_{1},...,%
\hat{y}_{j},...,y_{i-1}\right\rangle }(y_{i}-v)  \label{hoxij}
\end{equation}%
and $d_{n,t}\left( \hat{\jmath}\right) $\ be the Dickson algebra generator
of degree $2\left( p^{n-1}-p^{t}\right) $ in 
\begin{equation*}
\left( 
%TCIMACRO{\TeXButton{TeX field}{\mathbb F_p}}%
%BeginExpansion
\mathbb F_p%
%EndExpansion
[y_{1},\cdots ,\hat{y}_{j},...,y_{n}]\right) ^{GL(n-1,\mathbb{F}_{p})}
\end{equation*}

Let $\delta _{i,j}\in GL(n,\mathbb{F}_{p})$\ such that it permutes only the $%
i$ and $j$\ coordinates. Let\ 
\begin{equation}
h_{i}\left( j\right) :=\delta _{i,j}h_{i}=\prod\limits_{v\in \left\langle
y_{1},...,\hat{y}_{j},...,y_{i}\right\rangle }(y_{j}-v)  \label{h(i,j)}
\end{equation}%
for $j\leq i$.

\begin{lemma}
\label{h-hat}$h_{i}=h_{i}^{p}\left( \hat{\jmath}\right) -h_{i}\left( \hat{%
\jmath}\right) \left( h_{i-1}\left( j\right) \right) ^{p-1}$.
\end{lemma}

%TCIMACRO{\TeXButton{Proof}{\proof}}%
%BeginExpansion
\proof%
%EndExpansion
\begin{equation*}
h_{i}=\tprod\limits_{a}\prod\limits_{v\in \left\langle y_{2},\cdots
,y_{i-1}\right\rangle }(y_{i}-ay_{1}-v)=
\end{equation*}%
\begin{equation*}
\tprod\limits_{a}\tsum\limits_{t=0}^{i-2}\left( y_{i}+ay_{1}\right)
^{p^{i-2-t}}\left( -1\right) ^{t}d_{i-1,t}\left( \hat{1}\right)
=\tprod\limits_{a}\left( h_{i}\left( \hat{1}\right) +ah_{i-1}\left( 1\right)
\right)
\end{equation*}%
Since $\tsum\limits_{a}a\equiv 0\func{mod}p$, $\tsum\limits_{a_{i_{t}}\neq
a_{i_{l}}}\tprod\limits_{t=1}^{p-2}a_{i_{t}}\equiv 0\func{mod}p$\ and $%
\tprod\limits_{a\neq 0}a\equiv p-1\func{mod}p$, $h_{i}=h_{i}^{p}\left( \hat{1%
}\right) -h_{i}\left( \hat{1}\right) \left( h_{i-1}\left( 1\right) \right)
^{p-1}$. Now applying $\delta _{1,j}$ the statement follows.\ 
%TCIMACRO{\TeXButton{End Proof}{\endproof}}%
%BeginExpansion
\endproof%
%EndExpansion
\ 

The Dickson's result was extended for $H^{\ast }\left( V\right) ^{GL(2,%
\mathbb{F}_{p})}$ by Cardenas \ and Mui for the general case. For full
details please see \cite{Mui}.

In $E_{%
%TCIMACRO{\TeXButton{TeX field}{\mathbb Z}}%
%BeginExpansion
\mathbb Z%
%EndExpansion
}\left( x_{1},...,x_{n}\right) \otimes 
%TCIMACRO{\TeXButton{TeX field}{\mathbb Z}}%
%BeginExpansion
\mathbb Z%
%EndExpansion
[y_{1},...,y_{n}]$, let $M_{n,s_{1},...,s_{k}}$\ be defined as 
\begin{equation*}
\frac{1}{k!}\det \left( 
\begin{array}{cccccccccc}
x_{1} & \ldots & x_{1} & y_{1} & \ldots & \widehat{y}_{1}^{p^{s_{1}}} & 
\ldots & \widehat{y}_{1}^{p^{s_{k}}} & \ldots & y_{1}^{p^{n-1}} \\ 
\vdots &  & \vdots & \vdots &  & \vdots &  & \vdots &  & \vdots \\ 
x_{n} & \ldots & x_{n} & y_{n} & \ldots & \widehat{y}_{n}^{p^{s_{1}}} & 
\ldots & \widehat{y}_{n}^{p^{s_{k}}} & \ldots & y_{n}^{p^{n-1}}%
\end{array}%
\right)
\end{equation*}%
Here $0\leq s_{1}<...<s_{k}\leq n-1$. The columns $\left( 
\begin{array}{c}
\widehat{y}_{1}^{p^{s_{i}}} \\ 
\vdots \\ 
\widehat{y}_{n}^{p^{s_{i}}}%
\end{array}%
\right) $ are missing and the matrix for the proceeding determinant is filed
out with $k$\ columns of the form $\left( 
\begin{array}{c}
x_{1} \\ 
\vdots \\ 
x_{n}%
\end{array}%
\right) $\ to have $n$\ rows and columns. Let 
\begin{equation*}
M_{n,i}\left( \hat{t}\right) =Det\left( 
\begin{array}{cccccc}
x_{1} & y_{1} & \ldots & \widehat{y}_{1}^{p^{i}} & \ldots & y_{1}^{p^{n-2}}
\\ 
\ldots &  &  &  &  & \ldots \\ 
x_{n} & y_{n} & \ldots & \widehat{y}_{n}^{p^{i}} & \ldots & y_{n}^{p^{n-2}}%
\end{array}%
\right)
\end{equation*}%
and the $t$-th row is missing i.e. $\left[ x_{t},y_{t}...,y_{t}^{p^{n-2}}%
\right] $. Now the following formula is obvious:%
\begin{equation}
M_{n,n-1}=\left( -1\right) ^{t-1}[x_{t}L_{n,n-1}\left( \hat{t}\right)
-y_{t}M_{n,0}\left( \hat{t}\right) +...+\left( -1\right)
^{n-1}y_{t}^{p^{n-2}}M_{n,n-2}\left( \hat{t}\right) ]  \label{M(1)}
\end{equation}

We recall that $\hat{M}_{m,s_{1},...,s_{k}}=\omega M_{m,s_{1},...,s_{k}}$
and $\hat{d}_{m,t}=\omega d_{m,t}$ for $1\leq m\leq n$ and $\omega \in GL(n,%
\mathbb{F}_{p})$.

\begin{theorem}[Mui]
\label{mui} i) $H^{\ast }\left( V\right) ^{GL(n,\mathbb{F}_{p})}\cong
D_{n}\tbigoplus\limits_{k}\tbigoplus%
\limits_{s_{i}}D_{n}M_{n,s_{1},...,s_{k}}L_{n}^{p-2}$. Here a double
summation is taken over $k=1,...,n$\ and $0\leq s_{1}<...<s_{k}\leq n-1$.
Furthermore the generators satisfy: 1) $M_{n,s}^{2}=0$\ and \newline
2) $M_{n,s_{1}}...M_{n,s_{k}}=\left( -1\right)
^{k(k-1)/2}M_{n,s_{1},...,s_{k}}L_{n}^{k-1}$.

ii) $H^{\ast }\left( V\right) ^{U_{n}^{t}}\cong
H_{n}^{t}\tbigoplus\limits_{i}\tbigoplus\limits_{s_{t}}H_{n}^{t}\hat{M}%
_{i,s_{1},...,s_{k-1},i-1}$. Here $k\leq i\leq n$\ and $0\leq
s_{1}<...<s_{k-1}<i-1$.
\end{theorem}

The next lemma describes relations between exterior and polynomial algebra
generators.

\begin{lemma}
\label{relatio-M}i) Let $0\leq s_{1}<...<s_{k}\leq n-2$. Then 
\begin{equation*}
M_{n-1,s_{1},...,s_{k}}h_{n}=
\end{equation*}%
\begin{equation*}
M_{n,s_{1},...,s_{k}}-\tsum\limits_{\left( t_{1},...,t_{k}\right)
>(s_{k}-k+1,...,s_{k})}\left( -1\right) ^{k+i}M_{n,s_{1},...,\hat{s}%
_{i},...,s_{k}}d_{n-1,s_{i}}
\end{equation*}

ii) Let $0\leq s_{1}<...<s_{k}\leq k-1$. Then 
\begin{equation*}
M_{l,s_{1},...,s_{k}}h_{l+1}...h_{n}=M_{n,s_{1},...,s_{k}}+\tsum\limits_{%
\left( t_{1},...,t_{k}\right)
>(s_{k}-k+1,...,s_{k})}M_{n,t_{1},...,t_{k}}f_{t_{1},...,t_{k}}
\end{equation*}%
Here $f_{t_{1},...,t_{k}}\in H_{n}$.
\end{lemma}

The next Theorem is an extension of Mui's Theorem for parabolic subgroups (%
\cite{Kech1}).

\begin{theorem}[Kechagias]
\label{invar-Bock-Parab} Let $I=(n_{l},\cdots ,n_{1})$ be a sequence of
non-negative integers such that $\sum n_{i}=n$ and $P(I)$ be the associated
parabolic subgroup of $GL(n,\mathbb{F}_{p})$, then \ \ \ \ 
\begin{equation*}
H^{\ast }\left( V\right) ^{P\left( I\right) }\cong 
%TCIMACRO{\TeXButton{TeX field}{\mathbb F_p}}%
%BeginExpansion
\mathbb F_p%
%EndExpansion
(I)\tbigoplus\limits_{i}\tbigoplus\limits_{k}\tbigoplus\limits_{s_{t}}%
%TCIMACRO{\TeXButton{TeX field}{\mathbb F_p}}%
%BeginExpansion
\mathbb F_p%
%EndExpansion
(I)M_{\nu _{i},s_{1},...,s_{k}}L_{\nu _{i}}^{p-2}
\end{equation*}%
Here $1\leq i\leq \ell $, $\nu _{i}=\sum\limits_{t=1}^{i}n_{t}$, $1\leq
k\leq \nu _{i}$, $\nu _{i-1}\leq s_{k}$\ and \ \ $0\leq s_{1}<...<s_{k}\leq
\nu _{i}-1$.\ 
\end{theorem}

\section{The restriction map}

We remind the reader about a well known analogy between 
\begin{equation*}
U_{n}\leq B_{n}\leq P(I)\leq GL(n,\mathbb{F}_{p})
\end{equation*}%
and subgroups of the symmetric group $\Sigma _{p^{n}}$. There exists a
regular embedding \ $V\hookrightarrow \Sigma _{p^{n}}$\ which takes $u\in V$%
\ to the permutation on $V$\ induced by $v\mapsto u+v$.\ \ \ 

Let us recall that the wreath product between $H\leq \Sigma _{l}$ and $K\leq
\Sigma _{m}$ is defined by 
\begin{equation*}
1\rightarrow H^{m}\rightarrow K\int H\rightarrow K\rightarrow 1
\end{equation*}%
and $K\int H\leq \Sigma _{ml}$.

Let $\Sigma _{p^{n},p}:=\left( 
%TCIMACRO{\TeXButton{TeX field}{\mathbb Z_p}}%
%BeginExpansion
\mathbb Z_p%
%EndExpansion
\right) _{n}\idotsint \left( 
%TCIMACRO{\TeXButton{TeX field}{\mathbb Z_p}}%
%BeginExpansion
\mathbb Z_p%
%EndExpansion
\right) _{1}$ and $\Sigma (I):=\Sigma _{p^{n_{l}}}\int ...\int \Sigma
_{p^{n_{1}}}$. Then $\Sigma _{p^{n},p}$ is a $p$-Sylow subgroup of $\Sigma
_{p^{n}}$ and here is the analogy 
\begin{equation*}
\Sigma _{p^{n},p}\leq \Sigma (1,...,1)\leq \Sigma (I)\leq \Sigma _{p^{n}}
\end{equation*}%
Here the inclusion $V\hookrightarrow \Sigma _{p^{n},p}$ factors as follows 
\begin{equation*}
V=%
%TCIMACRO{\TeXButton{TeX field}{\mathbb Z_p}}%
%BeginExpansion
\mathbb Z_p%
%EndExpansion
\times \left( 
%TCIMACRO{\TeXButton{TeX field}{\mathbb Z_p}}%
%BeginExpansion
\mathbb Z_p%
%EndExpansion
\right) ^{n-1}\overset{1\times \Delta ^{p}}{\rightarrow }%
%TCIMACRO{\TeXButton{TeX field}{\mathbb Z_p}}%
%BeginExpansion
\mathbb Z_p%
%EndExpansion
\smallint \Sigma _{p^{n-1},p}\rightarrow \Sigma _{p}\smallint \Sigma
_{p^{n-1},p}\rightarrow \Sigma _{p^{n}}
\end{equation*}%
Moreover, the Weyl subgroups of $V$ in $\Sigma _{p^{n},p}$, $\Sigma (I)$,
and $\Sigma _{p^{n}}$ are the upper triangular group $U_{n}$, $P(I)$ and the
general linear group $GL(n,\mathbb{F}_{p})$ respectively. Please see \cite%
{Kuhn} Theorem 3.2.

Finally, $Aut(V)\cong GL(n,\mathbb{F}_{p})$ and let 
\begin{equation*}
\rho :W_{\Sigma _{p^{n}}}(V)\hookrightarrow GL(n,\mathbb{F}_{p})
\end{equation*}%
be the regular representation. Now the contragredient representation $\rho
^{\ast }$\ acts on $V^{\ast }\cong H^{1}\left( V\right) $. Here $\rho ^{\ast
}\left( g\right) =\rho \left( g^{-1}\right) ^{t}$. \ Moreover the Weyl
group, $W_{\Sigma _{p^{n}}}(V)$ $\cong GL(n,\mathbb{F}_{p})$, acts on $%
V^{\ast }$ as follows: 
\begin{equation*}
(a_{i,j})x_{k}:=\sum_{i}a_{i,k}x_{i}
\end{equation*}%
Here, $V^{\ast }=\left\langle x_{1},\cdots x_{n}\right\rangle $.

Let $E_{G}$\ and $B_{G}$\ denote the total and classifying spaces of a
finite group $G$. Let $H\leq G$\ be a subgroup, then $E_{G}$\ can also be a
total space for $H$\ and $pt\times _{H}E_{G}$\ is a model for $B_{G}$.
Moreover, 
\begin{equation*}
G/H\rightarrow B_{H}\overset{\pi }{\rightarrow }B_{G}
\end{equation*}%
\ is a fibration. The inclusion described above, $V\hookrightarrow G$,
induces a map \ 
\begin{equation*}
\left( res_{V}^{G}\right) ^{\ast }:H^{\ast }\left( G\right) \rightarrow
H^{\ast }\left( V\right) ^{W_{G}\left( V\right) }
\end{equation*}%
Here $G=\Sigma (I)$ and $H^{\ast }\left( G\right) :=H^{\ast }\left( B_{G},%
%TCIMACRO{\TeXButton{TeX field}{\mathbb Z}}%
%BeginExpansion
\mathbb Z%
%EndExpansion
/p\right) $. \ 

Since $H^{1}\left( V\right) \cong V^{\ast }$ and the Bockstein homomorphism
is an isomorphism $\beta :H^{1}\left( V\right) \rightarrow H^{2}\left(
V\right) $, let $y_{i}=\beta x_{i}$\ for \thinspace $1\leq i\leq n$. Now 
\begin{equation*}
H^{\ast }\left( V\right) =E_{%
%TCIMACRO{\TeXButton{TeX field}{\mathbb F_p}}%
%BeginExpansion
\mathbb F_p%
%EndExpansion
}\left( x_{1},\cdots x_{n}\right) \otimes 
%TCIMACRO{\TeXButton{TeX field}{\mathbb F_p}}%
%BeginExpansion
\mathbb F_p%
%EndExpansion
[y_{1},\cdots ,y_{n}]
\end{equation*}
and $H^{\ast }\left( V\right) ^{GL(n,\mathbb{F}_{p})}$ denotes the Dickson
algebra.

\begin{remark}
\label{transpose} Note that 
\begin{equation*}
\func{Im}\left( res_{V}^{G}\right) ^{\ast }\leq H^{\ast }\left( V\right)
^{W_{G}\left( V\right) }=\left( E_{%
%TCIMACRO{\TeXButton{TeX field}{\mathbb F_p}}%
%BeginExpansion
\mathbb F_p%
%EndExpansion
}\left( x_{1},\cdots x_{n}\right) \otimes 
%TCIMACRO{\TeXButton{TeX field}{\mathbb F_p}}%
%BeginExpansion
\mathbb F_p%
%EndExpansion
[y_{1},\cdots ,y_{n}]\right) ^{W_{G}\left( V\right) ^{t}}
\end{equation*}%
In other words we consider the transposes of the groups described above.
\end{remark}

The following important Theorem first proved by Cardenas for $n=2$\ and
extended by Kuhn provides the effective tools for our calculations. Here we
use a particular version of that Theorem. Please see VI, 1.6 in \cite{Ad-Mil}%
.

\begin{theorem}
(Cardenas, Mui, Kuhn).\newline
i) Let $res^{\ast }:H^{\ast }\left( \Sigma _{p}\smallint \Sigma
_{p^{n-1}}\right) \rightarrow H^{\ast }\left( V\right) $, then 
\begin{equation*}
\func{Im}\left( res^{\ast }\right) =H^{\ast }\left( V\right) ^{P(1,n-1)}\cap 
\func{Im}\left( res^{\ast }:H^{\ast }\left( \Sigma _{p^{n},p}\right)
\rightarrow H^{\ast }\left( V\right) \right)
\end{equation*}

ii) Let $res^{\ast }:H^{\ast }\left( \Sigma _{p^{n-1}}\smallint \Sigma
_{p}\right) \rightarrow H^{\ast }\left( V\right) $, then 
\begin{equation*}
\func{Im}\left( res^{\ast }\right) =H^{\ast }\left( V\right) ^{P(n-1,1)}\cap 
\func{Im}\left( res^{\ast }:H^{\ast }\left( \Sigma _{p^{n},p}\right)
\rightarrow H^{\ast }\left( V\right) \right)
\end{equation*}
\end{theorem}

Our first task is to give an invariant theoretic description of \newline
$\func{Im}\left( res^{\ast }:H^{\ast }\left( \Sigma _{p^{n},p}\right)
\rightarrow H^{\ast }\left( V\right) \right) $. Using a Theorem of Steenrod
and the action of the Steenrod algebra on upper triangular invariants we
compute this ring. For completeness we repeat some well known facts on group
cohomology. For full details please see VII in \cite{Steen}.

Let $H\lhd G$, then we have a fibering. An application of this fibering is
the following:%
\begin{equation*}
\left( B_{G}\right) ^{p}\overset{j}{\rightarrow }E_{%
%TCIMACRO{\TeXButton{TeX field}{\mathbb Z_p}}%
%BeginExpansion
\mathbb Z_p%
%EndExpansion
}\times _{%
%TCIMACRO{\TeXButton{TeX field}{\mathbb Z_p}}%
%BeginExpansion
\mathbb Z_p%
%EndExpansion
}\left( B_{G}\right) ^{p}\overset{\pi }{\rightarrow }B_{%
%TCIMACRO{\TeXButton{TeX field}{\mathbb Z_p}}%
%BeginExpansion
\mathbb Z_p%
%EndExpansion
}
\end{equation*}%
\ Here $G^{p}\lhd 
%TCIMACRO{\TeXButton{TeX field}{\mathbb Z_p}}%
%BeginExpansion
\mathbb Z_p%
%EndExpansion
\smallint G$ and $\left( B_{G}\right) ^{p}\simeq B_{G^{p}}$. The last
implies 
\begin{equation*}
H^{\ast }\left( G\right) \otimes ...\otimes H^{\ast }\left( G\right) \cong
H^{\ast }\left( G^{p}\right)
\end{equation*}

Let $\Delta ^{p}:B_{G}\rightarrow \left( B_{G}\right) ^{p}$\ be the diagonal
and 
\begin{equation*}
1\times \Delta ^{p}:B_{%
%TCIMACRO{\TeXButton{TeX field}{\mathbb Z_p}}%
%BeginExpansion
\mathbb Z_p%
%EndExpansion
}\times B_{G}\rightarrow B\left( 
%TCIMACRO{\TeXButton{TeX field}{\mathbb Z_p}}%
%BeginExpansion
\mathbb Z_p%
%EndExpansion
\smallint G\right) \simeq E_{%
%TCIMACRO{\TeXButton{TeX field}{\mathbb Z_p}}%
%BeginExpansion
\mathbb Z_p%
%EndExpansion
}\times _{%
%TCIMACRO{\TeXButton{TeX field}{\mathbb Z_p}}%
%BeginExpansion
\mathbb Z_p%
%EndExpansion
}\left( B_{G}\right) ^{p}
\end{equation*}%
the induced map. The image of the restriction map is the image of $1\times
\Delta ^{p}$. Now $H^{\ast }\left( 
%TCIMACRO{\TeXButton{TeX field}{\mathbb Z_p}}%
%BeginExpansion
\mathbb Z_p%
%EndExpansion
\smallint G\right) $ is an $H^{\ast }\left( 
%TCIMACRO{\TeXButton{TeX field}{\mathbb Z_p}}%
%BeginExpansion
\mathbb Z_p%
%EndExpansion
\right) $-module and $\left( \Delta ^{p}\right) ^{\ast }$ is an $H^{\ast
}\left( 
%TCIMACRO{\TeXButton{TeX field}{\mathbb Z_p}}%
%BeginExpansion
\mathbb Z_p%
%EndExpansion
\right) $-module homomorphism. Moreover the map $\pi ^{\ast }$ is a
monomorphism.

Let $\left\{ u_{j}|j\in J\right\} $ be an 
%TCIMACRO{\TeXButton{TeX field}{$\mathbb F_p$} }%
%BeginExpansion
$\mathbb F_p$
%EndExpansion
basis of $H^{\ast }\left( G\right) $. Then 
\begin{equation*}
M:=\left\langle u_{j}\otimes .....\otimes u_{j}|j\in J\right\rangle
\end{equation*}%
\ is an 
%TCIMACRO{\TeXButton{TeX field}{$\mathbb F_p$}}%
%BeginExpansion
$\mathbb F_p$%
%EndExpansion
-submodule of $H^{\ast }\left( G^{p}\right) $\ and 
\begin{equation*}
F:=\left\langle u_{j_{1}}\otimes ...\otimes u_{j_{p}}|j_{1}\leq ...\leq
j_{p}\ j_{1}<j_{p}\right\rangle
\end{equation*}%
\ is a free 
%TCIMACRO{\TeXButton{TeX field}{$\mathbb F_p$}}%
%BeginExpansion
$\mathbb F_p$%
%EndExpansion
-submodule of $H^{\ast }\left( G^{p}\right) $. It is well known that \ 
\begin{equation*}
H^{\ast }\left( 
%TCIMACRO{\TeXButton{TeX field}{\mathbb Z_p}}%
%BeginExpansion
\mathbb Z_p%
%EndExpansion
\smallint G\right) \cong H^{\ast }\left( 
%TCIMACRO{\TeXButton{TeX field}{\mathbb Z_p}}%
%BeginExpansion
\mathbb Z_p%
%EndExpansion
;\left( H^{\ast }\left( G\right) ^{p}\right) \right) \cong 
%TCIMACRO{\TeXButton{TeX field}{\mathbb F_p}}%
%BeginExpansion
\mathbb F_p%
%EndExpansion
\otimes F^{%
%TCIMACRO{\TeXButton{TeX field}{\mathbb Z_p}}%
%BeginExpansion
\mathbb Z_p%
%EndExpansion
}\oplus H^{\ast }\left( 
%TCIMACRO{\TeXButton{TeX field}{\mathbb Z_p}}%
%BeginExpansion
\mathbb Z_p%
%EndExpansion
\right) \otimes M
\end{equation*}%
Please see IV Theorem 1.7 in \cite{Ad-Mil}. If $\upsilon \in H^{\ast }\left( 
%TCIMACRO{\TeXButton{TeX field}{\mathbb Z_p}}%
%BeginExpansion
\mathbb Z_p%
%EndExpansion
\right) $, then $\upsilon $\ acts on $H^{\ast }\left( 
%TCIMACRO{\TeXButton{TeX field}{\mathbb Z_p}}%
%BeginExpansion
\mathbb Z_p%
%EndExpansion
\smallint G\right) $\ by $1^{p}\otimes \upsilon $.

Given a class $v\in H^{\ast }\left( G\right) $\ we have a class $v\otimes
.....\otimes v\in H^{\ast }\left( G^{p}\right) $. Now $\left( \Delta
^{p}\right) ^{\ast }\left( v\otimes ...\otimes v\right) =v^{p}$ and Steenrod
defined a map on the cochain level in order to compute the image of the
restriction map 
\begin{equation*}
P:H^{q}\left( G\right) \rightarrow H^{pq}\left( 
%TCIMACRO{\TeXButton{TeX field}{\mathbb Z_p}}%
%BeginExpansion
\mathbb Z_p%
%EndExpansion
\smallint G\right)
\end{equation*}%
such that \ $Pv$\ is the cohomology class $\varepsilon \otimes v^{p}$\ where 
$\varepsilon $ is the augmentation on the chain level. More precisely, 
\begin{equation*}
Pv=1\otimes v^{p}\in 
%TCIMACRO{\TeXButton{TeX field}{\mathbb F_p}}%
%BeginExpansion
\mathbb F_p%
%EndExpansion
\otimes F^{%
%TCIMACRO{\TeXButton{TeX field}{\mathbb Z_p}}%
%BeginExpansion
\mathbb Z_p%
%EndExpansion
}\oplus 
%TCIMACRO{\TeXButton{TeX field}{\mathbb F_p}}%
%BeginExpansion
\mathbb F_p%
%EndExpansion
\otimes M
\end{equation*}

Moreover, the Steenrod map satisfies 
\begin{equation*}
P\left( u\cup v\right) =\left( -1\right) ^{p(p-1)/2|u||v|}Pu\cup Pv
\end{equation*}%
Please see page 190 in \cite{Ad-Mil}. Now $H^{\ast }\left( 
%TCIMACRO{\TeXButton{TeX field}{\mathbb Z_p}}%
%BeginExpansion
\mathbb Z_p%
%EndExpansion
\right) \otimes \func{Im}P\cong H^{\ast }\left( 
%TCIMACRO{\TeXButton{TeX field}{\mathbb Z_p}}%
%BeginExpansion
\mathbb Z_p%
%EndExpansion
\right) \otimes M$\ \ and $H^{\ast }\left( 
%TCIMACRO{\TeXButton{TeX field}{\mathbb Z_p}}%
%BeginExpansion
\mathbb Z_p%
%EndExpansion
\right) \otimes \func{Im}\left( \Delta ^{p}\right) ^{\ast }P=\func{Im}\left(
\Delta ^{p}\right) ^{\ast }$.

\begin{theorem}
(Steenrod, May). Let $v\in H^{q}\left( G\right) $, $\eta =\left( p-1\right)
/2$\ and $\mu \left( q\right) =\left( \eta !\right) ^{-q}\left( -1\right)
^{\eta \left( q^{2}+q\right) /2}$. Then 
\begin{equation*}
\left( 1\times \Delta ^{p}\right) ^{\ast }Pv=\mu \left( q\right)
[\tsum\limits_{i}(-1)^{i}y^{(q-2i)\eta }\otimes
P^{i}v+\tsum\limits_{i}(-1)^{i+q}xy^{(q-2i)\eta -1}\otimes \beta P^{i}v]
\end{equation*}%
Here $H^{\ast }\left( 
%TCIMACRO{\TeXButton{TeX field}{\mathbb Z_p}}%
%BeginExpansion
\mathbb Z_p%
%EndExpansion
\right) \cong E_{%
%TCIMACRO{\TeXButton{TeX field}{\mathbb F_p}}%
%BeginExpansion
\mathbb F_p%
%EndExpansion
}\left( x\right) \otimes 
%TCIMACRO{\TeXButton{TeX field}{\mathbb F_p}}%
%BeginExpansion
\mathbb F_p%
%EndExpansion
[y]$.
\end{theorem}

Please see IV Theorem 4.1 in \cite{Ad-Mil}.

Now we are ready to prove the main Theorem of this section.

\begin{theorem}
\label{Steen-May} 
\begin{equation*}
\func{Im}\left( res^{\ast }:H^{\ast }\left( \Sigma _{p^{n},p}\right)
\rightarrow H^{\ast }\left( V\right) \right) \cong
\end{equation*}%
\begin{equation*}
E_{%
%TCIMACRO{\TeXButton{TeX field}{\mathbb F_p}}%
%BeginExpansion
\mathbb F_p%
%EndExpansion
}\left( \hat{M}_{1,0},\hat{M}_{2,1}\hat{L}_{1}^{\left( p-3\right) /2},\cdots
,\hat{M}_{n,n-1}\hat{L}_{n-1}^{\left( p-3\right) /2}\right) \otimes H_{n}^{t}
\end{equation*}
\end{theorem}

%TCIMACRO{\TeXButton{Proof}{\proof}}%
%BeginExpansion
\proof%
%EndExpansion
We apply induction on $n$. We shall prove\newline
i) $\left( 1\times \Delta ^{p}\right) ^{\ast }P\left( \hat{h}_{i}(\hat{n}%
)\right) =c\hat{h}_{i}$\ and\newline
ii) $\left( 1\times \Delta ^{p}\right) ^{\ast }P\left( \hat{M}_{i,i-1}(\hat{n%
})\hat{L}_{i-1}^{\left( p-3\right) /2}(\hat{n})\right) =c^{\prime }\hat{M}%
_{i,i-1}\hat{L}_{i-1}^{\left( p-3\right) /2}$.\newline
Here $c$, $c^{\prime }\in \left( 
%TCIMACRO{\TeXButton{TeX field}{\mathbb F_p}}%
%BeginExpansion
\mathbb F_p%
%EndExpansion
\right) ^{\ast }$. Or equivalently, \newline
$\left( 1\times \Delta ^{p}\right) ^{\ast }P\left( h_{i}(\hat{1})\right)
=ch_{i}\ $and \newline
$\left( 1\times \Delta ^{p}\right) ^{\ast }P\left( M_{i,i-1}(\hat{1}%
)L_{i-1}^{\left( p-3\right) /2}(\hat{1})\right) =c^{\prime
}M_{i,i-1}L_{i-1}^{\left( p-3\right) /2}$.\ \newline
i) We apply Steenrod-May's formula. 
\begin{equation}
\left( 1\times \Delta ^{p}\right) ^{\ast }P\left( h_{i}(\hat{1})\right) =\mu
\left( 2p^{i-2}\right) \tsum\limits_{m}(-1)^{m}y_{1}^{(2p^{i-2}-2m)\eta
}\otimes P^{m}h_{i}(\hat{1})  \label{h-Steen-May}
\end{equation}%
We recall definitions \ref{hoxij}, \ref{h(i,j)} and lemma \ref{h-hat}:%
\begin{equation}
h_{i}=h_{i}^{p}\left( \hat{1}\right) -h_{i}\left( \hat{1}\right) \left(
h_{i-1}\left( 1\right) \right) ^{p-1}  \label{hi}
\end{equation}%
The idea is to compare the coefficients of $y_{1}^{l}$ for certain $l$'s in
the expressions (\ref{h-Steen-May}) and (\ref{hi}).

We start with the action of Steenrod's algebra $P^{m}h_{i}(\hat{1})$. We
apply Theorem \ref{Steen-Action}\ repeatedly.

If $m=p^{i-2}$, then $P^{m}h_{i}(\hat{1})=h_{i}^{p}(\hat{1})$.

Now let $m=a_{i-3}p^{i-3}+...+a_{s}p^{s}$, then 
\begin{equation*}
P^{m}h_{i}(\hat{1})=\left( -1\right) ^{a_{i-3}}h_{i}(\hat{1})d_{i-1,i-2}(%
\hat{1})\binom{a_{i-3}+1}{a_{i-4}}\tprod\limits_{t=s}^{i-4}\binom{a_{t+1}}{%
a_{t}}\tprod\limits_{t=s}^{i-4}d_{i-1,t}^{a_{t}-a_{t-1}}(\hat{1})
\end{equation*}%
We recall definition \ref{h(i,j)}: 
\begin{equation*}
\left( h_{i-1}\left( 1\right) \right) ^{p-1}=\left( y_{1}^{p^{i-2}}+\tsum
\left( -1\right) ^{t}y_{1}^{p^{i-2-t}}d_{i-1,t}(\hat{1})\right) ^{p-1}
\end{equation*}%
Let $r\leq p-1$, $0\leq t_{1}<...<t_{r}\leq i-2$\ and $\lambda
_{t_{1}}+...+\lambda _{t_{r}}=p-1$. \ \ Then the coefficient of $%
y_{1}^{\Sigma \lambda _{t_{i}}p^{t_{i}}}$\ in the last expression is given
by \ \ \ 
\begin{equation*}
\left( -1\right) ^{\left( p-1\right) \left( i-2\right) -\Sigma \lambda
_{t_{i}}t_{i}}\QOVERD( ) {\left( p-1\right) !}{\lambda _{t_{1}}!...\lambda
_{t_{r}}!}\Pi d_{i-1,t_{i}}^{\lambda _{t_{i}}}(\hat{1})
\end{equation*}%
Here $\left( p-1\right) \left( i-2\right) -\Sigma \lambda
_{t_{i}}t_{i}\equiv \Sigma \lambda _{t_{i}}t_{i}\func{mod}2$.

Next the corresponding coefficient of $y_{1}$ in (\ref{h-Steen-May})\ shall
be considered.

Let $(p^{i-2}-m)(p-1)=\Sigma \lambda _{t_{i}}p^{t_{i}}$. Then 
\begin{equation*}
m(p-1)=p^{i-2}\left( p-1\right) -\Sigma \lambda _{t_{i}}p^{t_{i}}=\left(
b_{t_{r}}-1\right) p^{i-3}+b_{t_{r}}(p^{i-4}+...+p^{t_{r}})+
\end{equation*}%
\begin{equation*}
b_{t_{r-1}}(p^{t_{r}-1}+...+p^{t_{r-1}})+...+b_{t_{2}}(p^{t_{3}-1}+...+p^{t_{2}})+b_{t_{1}}(p^{t_{2}-1}+...+p^{t_{1}})
\end{equation*}%
Here $b_{t_{1}}=\lambda _{t_{1}}$, $b_{t_{2}}-b_{t_{1}}=\lambda _{t_{2}}$,
............., $b_{t_{r}}-b_{t_{r-1}}=\lambda _{t_{r}}$\ and $b_{t_{r}}=p-1$%
. Thus $b_{t_{i}}=a_{t_{i}}=...=a_{t_{i+1}-1}$ for $i\leq r-1$ and $%
b_{t_{r}}=a_{t_{r}}=...=a_{i-4}=a_{i-3}+1$. It is an easy computation to
prove that the exponents of $\left( -1\right) $ are equal in both sides i.e. 
$\left( a_{i-3}+m\right) \equiv \Sigma \lambda _{t_{i}}t_{i}\func{mod}2$.

We conclude $\left( 1\times \Delta ^{p}\right) ^{\ast }P\left( h_{i}(\hat{1}%
)\right) \equiv -\mu \left( 2p^{i-2}\right) h_{i}$.

ii) We shall prove that 
\begin{equation}
\left( 1\times \Delta ^{p}\right) ^{\ast }P\left( M_{i,i-1}(\hat{1}%
)L_{i-1}^{\left( p-3\right) /2}(\hat{1})\right) =c^{\prime
}M_{i,i-1}L_{i-1}^{\left( p-3\right) /2}  \label{Delt-M}
\end{equation}%
by comparing the corresponding coefficients of powers of $y_{1}$. First we
consider elements $\beta P^{m}\left( M_{i,i-1}(\hat{1})L_{i-1}^{\left(
p-3\right) /2}(\hat{1})\right) \neq 0$. Please see proposition \ref%
{Steen-ActionBock}. This\ is equivalent with \ \ \ 
\begin{equation}
m=p^{i-3}+...+1+\Sigma ^{l}a_{i_{t}}\left( p^{i-3}+...+p^{i_{t}}\right) 
\text{ and }\Sigma ^{l}a_{i_{t}}\leq \frac{p-3}{2}  \label{bhta-m}
\end{equation}%
In this case 
\begin{equation*}
\beta P^{m}\left( M_{i,i-1}(\hat{1})L_{i-1}^{\left( p-3\right) /2}(\hat{1}%
)\right) =\left( a_{i_{1}},...,a_{i_{l}}\right) L_{i}(\hat{1})\left( \tprod
L_{i-1,i_{t}}^{a_{t}}(\hat{1})\right) L_{i-1}^{a_{l+1}}
\end{equation*}%
Here $a_{l+1}=\left( \frac{p-3}{2}-\Sigma a_{i_{t}}\right) $ and 
\begin{equation*}
\left( a_{i_{1}},...,a_{i_{l}}\right) =\left( \left( p-3\right) /2\right)
!/\Sigma a_{i_{t}}!\left( \frac{p-3}{2}-\Sigma a_{i_{t}}\right) !
\end{equation*}%
In Steenrod-May's formula, the corresponding exponent of $y_{1}$ is 
\begin{equation*}
\frac{p-1}{2}\left( 2p^{i-2}-2\left( p^{i-3}+...+1\right) -2\left( \Sigma
a_{i_{t}}\left( p^{i-3}+...+p^{i_{t}}\right) \right) \right) -1=
\end{equation*}%
\begin{equation*}
\Sigma ^{l}a_{i_{t}}p^{i_{t}}+a_{l+1}p^{i-2}
\end{equation*}%
For each $m$\ satisfying condition (\ref{bhta-m}), \ 
\begin{equation*}
(-1)^{m+p^{i-2}}\beta P^{m}\left( M_{i,i-1}(\hat{1})L_{i-1}^{\left(
p-3\right) /2}(\hat{1})\right) x_{1}y_{1}^{\Sigma
a_{i_{t}}p^{i_{t}}+a_{l+1}p^{i-2}}=
\end{equation*}%
\begin{equation*}
\left( -1\right) ^{m+1}\left( a_{i_{1}},...,a_{i_{l}}\right) L_{i-1}(\hat{1}%
)\left( \tprod L_{i-2,i_{t}}^{a_{i_{t}}}(\hat{1})\right) L_{i-2}^{a_{l+1}}(%
\hat{1})x_{1}y_{1}^{\Sigma a_{i_{t}}p^{i_{t}}+a_{l+1}p^{i-2}}
\end{equation*}%
The corresponding coefficient of $x_{1}y_{1}^{\Sigma
a_{i_{t}}p^{i_{t}}+a_{l+1}p^{i-2}}$ in the decomposition of\newline
$M_{i,i-1}L_{i-1}^{\left( p-3\right) /2}$\ (right hand side in (\ref{Delt-M}%
)) with respect to $x_{1}y_{1}$\ (according to formulas \ref{L(1)} and \ref%
{M(1)}) is 
\begin{equation*}
\left( -1\right) ^{\Sigma _{t=1}^{l+1}a_{i_{t}}i_{t}}\left(
a_{i_{1}},...,a_{i_{l}}\right) L_{i-1}(\hat{1})\left( \tprod
L_{i-1,i_{t}}^{a_{i_{t}}}(\hat{1})\right) L_{i-1,i-2}^{a_{l+1}}(\hat{1})
\end{equation*}%
Those two elements differ by $\left( -1\right) ^{1+\left( i-2\right) \left(
p-1\right) /2}$. \ 

Next we consider elements of the form 
\begin{equation*}
\left( -1\right) ^{m}P^{m}\left( M_{i,i-1}(\hat{1})L_{i-1}^{\left(
p-3\right) /2}(\hat{1})\right) y_{1}^{\frac{p-1}{2}\left( 2p^{i-2}-2m\right)
}
\end{equation*}%
in the left hand side of (\ref{Delt-M}).

For non-zero elements we have $m=p^{i-3}+...+p^{k}+m^{\prime }$ with \newline
$m^{\prime }=\Sigma ^{l}a_{i_{t}}\left( p^{i-3}+...+p^{i_{t}}\right) $ and $%
\Sigma ^{l}a_{i_{t}}\leq \frac{p-3}{2}$. Replacing $m$\ in the exponent of $%
y_{1}$\ it takes the form 
\begin{equation*}
\frac{p-1}{2}\left( 2p^{i-2}-2\left( p^{i-3}+...+1\right) -2\left( \Sigma
a_{i_{t}}\left( p^{i-3}+...+p^{i_{t}}\right) \right) \right) +p^{k}-1
\end{equation*}%
As before the corresponding coefficients of $y_{1}$ to the particular
exponent differ by 
\begin{equation*}
\left( -1\right) ^{1+\left( i-2\right) \left( p-1\right) /2}:m-\left(
k-1+\Sigma ^{l+1}a_{i_{t}}i_{t}\right) \equiv 1+\left( i-2\right) \left(
p-1\right) /2\func{mod}p
\end{equation*}%
Now the proof is complete.\ \ \ 
%TCIMACRO{\TeXButton{End Proof}{\endproof}}%
%BeginExpansion
\endproof%
%EndExpansion

\begin{proposition}
\label{decomp-restric} The image $\func{Im}\left( res^{\ast }:H^{\ast
}\left( \Sigma _{p^{n},p}\right) \rightarrow H^{\ast }\left( V\right)
\right) $ is isomorphic with 
\begin{equation*}
H_{n}^{t}\tbigoplus\limits_{i}\tbigoplus\limits_{s_{t}}H_{n}^{t}\hat{M}%
_{i,s_{1},...,s_{k-1},i-1}\hat{L}_{i-1}^{\left( p-3\right)
/2}\tprod\limits_{1}^{k-1}\hat{L}_{s_{t}}^{\left( p-3\right)
/2}\tprod\limits_{1}^{k-1}\hat{L}_{s_{t}+1}^{\left( p-3\right) /2}
\end{equation*}%
Here $k\leq i\leq n$\ and $0\leq s_{1}<...<s_{k-1}<i-1$.
\end{proposition}

%TCIMACRO{\TeXButton{Proof}{\proof}}%
%BeginExpansion
\proof%
%EndExpansion
This is an application of Theorem \ref{Steen-May} and lemma \ref{relatio-M}.%
%TCIMACRO{\TeXButton{End Proof}{\endproof}}%
%BeginExpansion
\endproof%
%EndExpansion

The next Theorem is an application of last Theorem and Cardenas-Mui-Kuhn
Theorem.

\begin{theorem}
\label{Restr-Parab}$\func{Im}\left( res^{\ast }:H^{\ast }\left( \Sigma
_{p^{n_{l}}}\smallint ...\smallint \Sigma _{p^{n_{1}}}\right) \rightarrow
H^{\ast }\left( V\right) \right) $ is isomorphic to the subalgebra generated
by 
\begin{equation*}
\left\{ 
\begin{array}{c}
\hat{d}_{\nu _{i},\nu _{i}-k_{i}}\ ,\hat{M}_{\nu _{i},\nu _{i}-k_{i}}\left( 
\hat{L}_{\nu _{i}}\right) ^{p-2},\hat{M}_{\nu _{i},\nu _{i}-k_{j},\nu
_{i}-k_{i}}\left( \hat{L}_{\nu _{i}}\right) ^{p-2}|\  \\ 
1\leq i\leq \ell ,\ 1\leq k_{i}\leq n_{i},k_{i}<k_{j}<\nu _{i},\ \nu
_{i}=\sum\limits_{t=1}^{i}n_{t}%
\end{array}%
\right\}
\end{equation*}%
Subject to relations described in Theorem \ref{mui}\ and lemma \ref%
{relatio-M}.\ 
\end{theorem}

\section{Relations between parabolic and Dickson algebra generators}

Since $D_{n}$ is a subalgebra of 
%TCIMACRO{\TeXButton{TeX field}{$\mathbb F_p$}}%
%BeginExpansion
$\mathbb F_p$%
%EndExpansion
$[V]^{P(I)}$, any Dickson algebra generator can be decomposed in terms of
generators of the later algebra. We shall describe these relations in this
section for $I=(n-1,1)$ and $\left( 1,n-1\right) $.

We recall that a Dickson algebra generator $d_{n,n-i}$ consists of the sum
of all possible combinations of $i$ elements from $\{h_{1}^{p-1},\cdots
,h_{n}^{p-1}\}$ in certain p-th exponents (proposition \ref{formula}) and
this might be more than what a $P(I)$-generator needs. For instance, we
would like to replace $d_{n,n-1}$ by another element which is a $P(I)$%
-invariant but not a $GL(n,\mathbb{F}_{p})$-one. An example is in order.

\textbf{Example} \label{Ex1} Proposition \ref{formula} is applied.

a) Let $n=4$ and $n_{1}=3$. 
\begin{equation*}
d_{4,3}=h_{1}^{(p-1)p^{3}}+h_{2}^{(p-1)p^{2}}+h_{3}^{(p-1)p}+h_{4}^{(p-1)}%
\Rightarrow
\end{equation*}%
\begin{equation*}
d_{4,3}-h_{1}^{(p-1)p^{3}}-h_{2}^{(p-1)p^{2}}-h_{3}^{(p-1)p}=h_{4}^{(p-1)}
\end{equation*}%
And this polynomial is a $P(3,1)$-invariant.

b) Let $n=4$ and $n_{1}=1$. 
\begin{equation*}
d_{4,3}=h_{1}^{(p-1)p^{3}}+h_{2}^{(p-1)p^{2}}+h_{3}^{(p-1)p}+h_{4}^{(p-1)}%
\Rightarrow
\end{equation*}%
\begin{equation*}
d_{4,3}-h_{1}^{(p-1)p^{3}}=h_{2}^{(p-1)p^{2}}+h_{3}^{(p-1)p}+h_{4}^{(p-1)}
\end{equation*}%
And this polynomial in a $P(1,3)$-invariant. Let us call the last sum $%
d_{4,3}\left( I\right) $. Thus 
\begin{equation*}
d_{4,3}\left( I\right) =d_{4,3}-h_{1}^{(p-1)p^{3}}
\end{equation*}%
Next we consider $d_{4,2}$: 
\begin{equation*}
d_{4,2}=h_{1}^{(p-1)p^{2}}h_{2}^{(p-1)p^{2}}+h_{1}^{(p-1)p}h_{3}^{(p-1)p^{2}}+h_{1}^{(p-1)}h_{4}^{(p-1)p^{2}}+h_{2}^{(p-1)p}h_{3}^{(p-1)p}+
\end{equation*}%
\begin{equation*}
h_{2}^{(p-1)}h_{4}^{(p-1)p}+h_{3}^{(p-1)}h_{4}^{(p-1)}\Rightarrow \ 
\end{equation*}%
\begin{equation*}
d_{4,2}-(h_{2}^{(p-1)p^{3}}+h_{3}^{(p-1)p}+h_{4}^{(p-1)})h_{1}^{(p-1)p^{2}}=h_{3}^{(p-1)}h_{4}^{(p-1)}
\end{equation*}%
Let us call the last sum $d_{4,2}\left( I\right) $. Thus 
\begin{equation*}
d_{4,2}\left( I\right) =d_{4,2}-h_{1}^{(p-1)p^{2}}d_{4,3}\left( I\right)
\end{equation*}%
Now $d_{4,1}$: 
\begin{equation*}
d_{4,1}=h_{2}^{(p-1)}h_{3}^{(p-1)}h_{4}^{(p-1)}+h_{1}^{(p-1)}h_{3}^{(p-1)}h_{4}^{(p-1)p}+h_{1}^{(p-1)p}h_{2}^{(p-1)p}h_{4}^{(p-1)p}+
\end{equation*}%
\begin{equation*}
h_{1}^{(p-1)p}h_{2}^{(p-1)p}h_{3}^{(p-1)p}\Rightarrow
\end{equation*}%
\begin{equation*}
d_{4,1}-(h_{3}^{(p-1)}h_{4}^{(p-1)p}+h_{2}^{(p-1)p}h_{4}^{(p-1)p}+h_{2}^{(p-1)p}h_{3}^{(p-1)p})h_{1}^{(p-1)p}=
\end{equation*}%
\begin{equation*}
h_{2}^{(p-1)}h_{3}^{(p-1)}h_{4}^{(p-1)}
\end{equation*}%
Thus 
\begin{equation*}
d_{4,1}\left( I\right) =d_{4,1}-h_{1}^{(p-1)p}d_{4,2}\left( I\right)
\end{equation*}%
Finally $d_{4,0}$: 
\begin{equation*}
d_{4,0}=h_{1}^{(p-1)}h_{2}^{(p-1)}h_{3}^{(p-1)}h_{4}^{(p-1)}=h_{1}^{(p-1)}d_{4,1}\left( I\right)
\end{equation*}

\begin{remark}
\label{ParabolicGener} According to proposition \ref{formula} each Dickson
algebra generator is a function on $\left\{
h_{1}^{p-1},...,h_{n}^{p-1}\right\} $: $%
F_{n,i}(h_{1}^{p-1},...,h_{n}^{p-1})=d_{n,n-i}$. Let $I=\left(
n_{1},n-n_{1}\right) $ we define 
\begin{equation*}
d_{n,n-i}\left( I\right) =F_{n,i}(0,...,0,h_{i+1}^{p-1},...,h_{n}^{p-1})
\end{equation*}%
for $n-n_{1}\geq i$.
\end{remark}

Let us note that $d_{n,n-i}\left( I\right) $ also depends on the value of $%
n_{1}$. Moreover, the new polynomial is a summand of $d_{n,n-i}$ and it will
be expressed in terms of old generators. The following proposition is an
application of corollary \ref{corol1} and Theorem \ref{invar-Bock-Parab}.

\begin{proposition}
Let $I=(n-1,1)$, then 
\begin{equation*}
%TCIMACRO{\TeXButton{TeX field}{\mathbb F_p}}%
%BeginExpansion
\mathbb F_p%
%EndExpansion
\left( I\right) =%
%TCIMACRO{\TeXButton{TeX field}{\mathbb F_p}}%
%BeginExpansion
\mathbb F_p%
%EndExpansion
[d_{n-1,i},h_{n}^{p-1}\ |\ 0\leq i\leq n-2]
\end{equation*}
and 
\begin{equation*}
H^{\ast }\left( V\right) ^{P\left( I\right) }\cong 
%TCIMACRO{\TeXButton{TeX field}{\mathbb F_p}}%
%BeginExpansion
\mathbb F_p%
%EndExpansion
\left( I\right) \tbigoplus\limits_{s_{i}}%
%TCIMACRO{\TeXButton{TeX field}{\mathbb F_p}}%
%BeginExpansion
\mathbb F_p%
%EndExpansion
\left( I\right) \left[ M_{n-1,s_{1},...,s_{k}}L_{n-1}^{p-2}\tbigoplus%
\limits_{t_{i}}M_{n,t_{1},...,t_{k},n-1}h_{n}^{p-2}\right]
\end{equation*}%
Here \ $0\leq s_{1}<...<s_{k}\leq n-2$ and $0\leq t_{1}<...<t_{k}<n-1$.
\end{proposition}

Next we compute the action of Steenrod's algebra on an upper triangular
generator.

\begin{theorem}
\label{Steen-Action}i) Let $m=\tsum\limits_{t=s}^{n-2}a_{t}p^{t}$ and $%
p-1\geq a_{n-2}+1\geq a_{n-3}\geq ...\geq a_{s}\geq 0$, then 
\begin{equation*}
P^{m}h_{n}=h_{n}\left( d_{n-1,n-2}\left( -1\right) ^{a_{n-2}}\binom{a_{n-2}+1%
}{a_{n-3}}\tprod\limits_{t=s}^{n-3}\binom{a_{t}}{a_{t-1}}\tprod%
\limits_{t=s}^{n-2}d_{n-1,t}^{a_{t}-a_{t-1}}\right)
\end{equation*}%
ii) Let $m=p^{n-1}$, then \ 
\begin{equation*}
P^{m}h_{n}=h_{n}^{p}
\end{equation*}%
iii) For all other cases, $P^{m}h_{n}=0$.
\end{theorem}

%TCIMACRO{\TeXButton{Proof}{\proof}}%
%BeginExpansion
\proof%
%EndExpansion
$P^{m}h_{n}=P^{m}\left( \tsum\limits_{t=0}^{n-1}\left( -1\right)
^{t}y_{n}^{p^{n-1-t}}d_{n-1,t}\right) =$

$\tsum\limits_{t=0}^{n-1}\left( -1\right)
^{t}\tsum\limits_{i+j=m}P^{i}y_{n}^{p^{n-1-t}}P^{j}d_{n-1,t}$.\ If $%
P^{i}y_{n}^{p^{n-1-t}}=0$\ for all $t$, then $P^{m}h_{n}=0$. Thus $%
P^{m}h_{n}\neq 0$\ implies $\exists t$ and $i$\ such that 
\begin{equation*}
P^{i}y_{n}^{p^{n-1-t}}=\left\{ 
\begin{array}{c}
y_{n}^{p^{n-1-t}}\text{ for }i=0 \\ 
P^{i}y_{n}^{p^{n-t}}\text{ for }i=p^{n-1-t} \\ 
0\text{, otherwise}%
\end{array}%
\right. \text{ and }P^{j}d_{n-1,t}\neq 0
\end{equation*}%
Thus $y_{n}$ divides $P^{m}h_{n}$.\ Since $P^{m}h_{n}\in H_{n}$, $%
P^{m}h_{n}=h_{n}f$ and $f\in H_{n}$.

Let $m<p^{n-1}$\ and $P^{m}h_{n}\neq 0$, then $P^{m}y_{n}^{p^{n-1}}=0$ and $%
P^{i}y_{n}^{p^{n-2}}\neq 0$. Otherwise, $P^{i}y_{n}^{p^{l}}=y_{n}^{p^{l+1}}$
for $l+1<n-1$. In that case $P^{m}h_{n}=0$.\ Thus $i=p^{n-2}$\ and $%
m=p^{n-2}+j$.

According to Theorem \ref{Steenr-action-dn}, $P^{j}d_{n-1,n-2}\neq 0$\ if
and only if $j=\tsum\limits_{t=s}^{n-2}a_{t}p^{t}$ and $p-1\geq
a_{n-2}+1\geq a_{n-3}\geq ...\geq a_{s}\geq 0$. Now the statement follows.%
%TCIMACRO{\TeXButton{End Proof}{\endproof}}%
%BeginExpansion
\endproof%
%EndExpansion
\ 

\begin{proposition}
\cite{Kech1}\label{Steen-ActionBock}i)$P^{m}M_{i,i-1}L_{i-1}^{\left(
p-3\right) /2}=$%
\begin{equation*}
M_{i,i-1}P^{m}L_{i-1}^{\left( p-3\right) /2}+\Sigma
_{0}^{i-1}M_{i,t}P^{m-\left( p^{i-1}+...+p^{t}\right) }L_{i-1}^{\left(
p-3\right) /2}
\end{equation*}%
\linebreak ii) If $m=p^{i-2}+...+1+\Sigma a_{i_{t}}(p^{i-2}+...+p^{i_{t}})$
and $\Sigma a_{i_{t}}\leq \frac{p-3}{2}$, then\ 

\begin{equation*}
\beta P^{m}M_{i,i-1}L_{i-1}^{\left( p-3\right) /2}=\left(
a_{i_{1}},...,a_{i_{l}}\right) L_{i}\left( \tprod
L_{i-1,i_{t}}^{a_{t}}\right) L_{i-1}^{\left( p-3\right) /2-\Sigma a_{i_{t}}}
\end{equation*}%
And $\beta P^{m}M_{i,i-1}L_{i-1}^{\left( p-3\right) /2}=0$, otherwise.

Here $\left( a_{i_{1}},...,a_{i_{l}}\right) =\left( \left( p-3\right)
/2\right) !/\Sigma a_{i_{t}}!\left( \frac{p-3}{2}-\Sigma a_{i_{t}}\right) !$.
\end{proposition}

According to Theorem \ref{Steen-Action}\ and proposition \ref%
{Steen-ActionBock},\ the action of Steenrod's algebra is closed on the
generating set above.

Next we proceed to the case $I=\left( 1,n-1\right) $.

\begin{proposition}
\label{decompos}Let $I=\left( 1,n-1\right) $ and $n-1\geq i\geq 1$, then $%
d_{n,n-i}\left( I\right) $ can be decomposed in terms of $d_{n,n-t}$ and
vise versa for $t<i$ as follows 
\begin{eqnarray*}
d_{n,n-i} &=&h_{1}^{(p-1)p^{n-i}}d_{n,n-i+1}\left( I\right) +d_{n,n-i}\left(
I\right) \\
d_{n,n-i}\left( I\right) &=&d_{n,n-i}-\sum\limits_{t=n-i}^{n-1}\left(
-1\right) ^{t+i+1-n}h_{1}^{(p-1)(p^{n-i}+...+p^{t})}d_{n,t+1}
\end{eqnarray*}
\end{proposition}

%TCIMACRO{\TeXButton{Proof}{\proof}}%
%BeginExpansion
\proof%
%EndExpansion
We apply proposition \ref{formula} and induction on $i$. $d_{n,n-i}$ is a
combination of $i$ elements from $\{h_{1}^{(p-1)},...,h_{n}^{(p-1)}\}$ or $1$
from $\{h_{1}^{(p-1)}\}$ and $i-1$ from $\{h_{2}^{(p-1)},...,h_{n}^{(p-1)}\}$
on certain powers: 
\begin{equation*}
d_{n,n-i}=h_{1}^{(p-1)p^{n-i}}d_{n,n-i+1}\left( I\right) +d_{n,n-i}\left(
I\right)
\end{equation*}%
Now the claim follows. 
%TCIMACRO{\TeXButton{End Proof}{\endproof}}%
%BeginExpansion
\endproof%
%EndExpansion

\begin{theorem}
\label{d'}Let $I=\left( 1,n-1\right) $, then 
\begin{equation*}
%TCIMACRO{\TeXButton{TeX field}{\mathbb F_p}}%
%BeginExpansion
\mathbb F_p%
%EndExpansion
\left( I\right) =%
%TCIMACRO{\TeXButton{TeX field}{\mathbb F_p}}%
%BeginExpansion
\mathbb F_p%
%EndExpansion
[h_{1}^{p-1},d_{n,i}\left( I\right) \ |\ 1\leq i\leq n-1]
\end{equation*}
and 
\begin{equation*}
H^{\ast }\left( V\right) ^{P\left( I\right) }\cong 
%TCIMACRO{\TeXButton{TeX field}{\mathbb F_p}}%
%BeginExpansion
\mathbb F_p%
%EndExpansion
\left( I\right) \oplus 
%TCIMACRO{\TeXButton{TeX field}{\mathbb F_p}}%
%BeginExpansion
\mathbb F_p%
%EndExpansion
\left( I\right) \left[ M_{1,0}h_{1}^{p-2}\tbigoplus%
\limits_{t_{l}}M_{n,t_{1},...,t_{k}}L_{n}^{p-2}\right]
\end{equation*}%
Here $1\leq t_{k}$ and\ $0\leq t_{1}<...<t_{k}\leq n-1$.
\end{theorem}

%TCIMACRO{\TeXButton{Proof}{\proof}}%
%BeginExpansion
\proof%
%EndExpansion
Because of last proposition the $d_{n,i}\left( I\right) $'s\ are invariants
and consist a polynomial basis. The claim follows from Theorem \ref%
{invar-Bock-Parab}.\ 
%TCIMACRO{\TeXButton{End Proof}{\endproof}}%
%BeginExpansion
\endproof%
%EndExpansion

\begin{proposition}
The action of Steenrod's algebra on the generating set \newline
$\left\{ h_{1}^{p-1},d_{n,i}\left( I\right) \ |\ 1\leq i\leq n-1\right\} $
is closed.
\end{proposition}

\begin{proof}
We need to evaluate the action of Steenrod's algebra for the Steenrod
algebra generators $P^{p^{l}}$\ only. We recall Theorem \ref%
{Steenr-action-dn}.\ \ 
\begin{equation*}
P^{p^{l}}d_{n,i}=\left\{ 
\begin{array}{c}
d_{n,i-1}\text{, for }l=i-1 \\ 
-d_{n,i}d_{n,n-1}\text{, for }l=n-1 \\ 
0\text{, otherwise}%
\end{array}%
\right. \text{ and }
\end{equation*}%
\begin{equation*}
P^{p^{l}}h_{1}^{\left( p-1\right) p^{k}}=\left\{ 
\begin{array}{c}
-h_{1}^{p^{k+1}-\left( p-2\right) p^{k}}\text{, for }l=k \\ 
0\text{, otherwise}%
\end{array}%
\right.
\end{equation*}%
Because of proposition \ref{decompos}, we have to consider $%
P^{p^{n-1}},...,P^{p^{n-i-1}}$ only. Let $n-i>1$, then $P^{p^{l}}d_{n,i}%
\left( I\right) $ is a function on the set: 
\begin{equation*}
\left\{ h_{1}^{p-1},d_{n,i}\left( I\right) \ |\ 1\leq i\leq n-1\right\}
\end{equation*}%
Let $n-i=1$, then\ we apply relation $d_{n,0}=h_{1}^{(p-1)}d_{n,1}\left(
I\right) $ on $P^{p^{l}}d_{n,1}\left( I\right) $.
\end{proof}

For the general case please see \cite{Kech3}.

\section{$%
%TCIMACRO{\TeXButton{TeX field}{\mathbb F_p}}%
%BeginExpansion
\mathbb F_p%
%EndExpansion
(n-1,1)$ and $%
%TCIMACRO{\TeXButton{TeX field}{\mathbb F_p}}%
%BeginExpansion
\mathbb F_p%
%EndExpansion
(1,n-1)$ as free modules over $D_{n}$}

$D_{n}$ serves as a homogeneous system of parameters and in fact both 
%TCIMACRO{\TeXButton{TeX field}{$\mathbb F_p$}}%
%BeginExpansion
$\mathbb F_p$%
%EndExpansion
$[V]^{U_{n}}$ and 
%TCIMACRO{\TeXButton{TeX field}{$\mathbb F_p$}}%
%BeginExpansion
$\mathbb F_p$%
%EndExpansion
$(I)$ are free $D_{n}$-modules. A free basis has been given for 
%TCIMACRO{\TeXButton{TeX field}{$\mathbb F_p$}}%
%BeginExpansion
$\mathbb F_p$%
%EndExpansion
$[V]^{U_{n}}$ as a module over $D_{n}$ ( \cite{C-H} and \cite{Kech1}).

Since $U_{n}$ is a $p$-Sylow subgroup of $GL(n,\mathbb{F}_{p})$ and $H_{n}$
is a polynomial algebra, 
%TCIMACRO{\TeXButton{TeX field}{$\mathbb F_p$}}%
%BeginExpansion
$\mathbb F_p$%
%EndExpansion
$(I)$ is Cohen-Macaulay. Hence, 
%TCIMACRO{\TeXButton{TeX field}{$\mathbb F_p$}}%
%BeginExpansion
$\mathbb F_p$%
%EndExpansion
$(I)$ is a free module over $D_{n}$.

\begin{remark}
\label{FreeModule}i) The rank of 
%TCIMACRO{\TeXButton{TeX field}{$\mathbb F_p$}}%
%BeginExpansion
$\mathbb F_p$%
%EndExpansion
$(I)$ over $D_{n}$ is $\left[ GL(n,\mathbb{F}_{p}):P(I)\right] $. \newline
ii) Let $P(G,t)$ denote the Poincar\'{e} series of 
%TCIMACRO{\TeXButton{TeX field}{$\mathbb F_p$}}%
%BeginExpansion
$\mathbb F_p$%
%EndExpansion
$(n-1,1)$. Note that $|d_{n-1,i}|=p^{n-1}-p^{i}\,$ divides $|d_{n-1+1,i+1}|$
and hence \newline
$P\left( D_{n},t\right) /\,P\left( G,t\right) =$ $\prod\limits_{i}\left(
1+t^{|d_{n-1,i}|}+t^{2|d_{n-1,i}|}+\cdots +t^{\left( p-1\right)
|d_{n-1,i}|}\right) $.
\end{remark}

\begin{definition}
Let the symbol $B_{A}\left( A^{\prime }\right) $ stand for a free module
basis of the algebra $A^{\prime }$ over the algebra $A$.
\end{definition}

\begin{theorem}
\cite{C-H} $B_{D_{n}}\left( H_{n}\right) =\left\{
h_{1}^{r_{1}}...h_{n}^{r_{n}}\ |\ 0\leq r_{i}<p^{n-i+1}-1\right\} $ is a
free module basis for $H_{n}$\ over $D_{n}$.\ 
\end{theorem}

\begin{corollary}
\label{Corol-Restr-Bock-Upper}$\func{Im}\left( res^{\ast }:H^{\ast }\left(
\Sigma _{p^{n},p}\right) \rightarrow H^{\ast }\left( V\right) \right) $ is
isomorphic to the free module over $D_{n}$ on \ 
\begin{equation*}
\left\{ 
\begin{array}{c}
\hat{M}_{i,s_{1},...,s_{k-1},i-1}\hat{L}_{i-1}^{\left( p-3\right)
/2}\tprod\limits_{2}^{k}\hat{L}_{s_{t}}^{\left( p-3\right)
/2}\tprod\limits_{2}^{k}\hat{L}_{s_{t}+1}^{\left( p-3\right) /2}\hat{h}%
_{1}^{r_{1}}...\hat{h}_{n}^{r_{n}}\ |\  \\ 
0\leq r_{i}<p^{n-i+1}-1,k\leq i\leq n,0\leq s_{1}<...<s_{k-1}<i-1%
\end{array}%
\right\}
\end{equation*}
\end{corollary}

%TCIMACRO{\TeXButton{Proof}{\proof}}%
%BeginExpansion
\proof%
%EndExpansion
This is an application of last Theorem and proposition \ref{decomp-restric}.%
%TCIMACRO{\TeXButton{End Proof}{\endproof}}%
%BeginExpansion
\endproof%
%EndExpansion

\begin{proposition}
\label{1,n-1}$B_{D_{n}}\left( 
%TCIMACRO{\TeXButton{TeX field}{\mathbb F_p}}%
%BeginExpansion
\mathbb F_p%
%EndExpansion
(1,n-1)\right) =\left\{ h_{1}^{(p-1)m}\,|\,0\leq m\leq A_{1}\right\} $ is a
free module basis for $%
%TCIMACRO{\TeXButton{TeX field}{\mathbb F_p}}%
%BeginExpansion
\mathbb F_p%
%EndExpansion
(1,n-1)$ over $D_{n}$. Here $A_{1}=p^{n-1}+...+p$.
\end{proposition}

%TCIMACRO{\TeXButton{Proof}{\proof}}%
%BeginExpansion
\proof%
%EndExpansion
Our statement follows directly from the following formulas: 
\begin{equation*}
d_{n,0}=d_{n,1}\left( I\right) h_{1}^{(p-1)}
\end{equation*}%
\begin{equation*}
d_{n,1}=d_{n,1}\left( I\right)
+\sum\limits_{t=1}^{n-1}(-1)^{t}d_{n,1+t}h_{1}^{(p-1)p^{t}+...+p}
\end{equation*}%
\begin{equation*}
d_{n,0}=d_{n,1}h_{1}^{(p-1)}+\sum%
\limits_{t=1}^{n-1}(-1)^{t}d_{n,1+t}h_{1}^{(p-1)p^{t}+p^{t-1}+...+1}
\end{equation*}%
%TCIMACRO{\TeXButton{End Proof}{\endproof}}%
%BeginExpansion
\endproof%
%EndExpansion

\begin{corollary}
\cite{Kech1} 
\begin{equation*}
B_{D_{n}}\left( 
%TCIMACRO{\TeXButton{TeX field}{\mathbb F_p}}%
%BeginExpansion
\mathbb F_p%
%EndExpansion
(1,...,1)\right) =\left\{
h_{1}^{(p-1)m_{1}}...h_{n-1}^{(p-1)m_{n-1}}\,|\,0\leq m_{i}\leq A_{i}\right\}
\end{equation*}
is a free module basis for $%
%TCIMACRO{\TeXButton{TeX field}{\mathbb F_p}}%
%BeginExpansion
\mathbb F_p%
%EndExpansion
(1,...,1)$ over $D_{n}$. Here $A_{i}=p^{n-i}+...+p$.
\end{corollary}

In the opposite direction as in the last proposition, we consider the
analogue statement. Next lemma demonstrates our approach.

\begin{lemma}
$B_{D_{4}}\left( 
%TCIMACRO{\TeXButton{TeX field}{\mathbb F_p}}%
%BeginExpansion
\mathbb F_p%
%EndExpansion
(3,1)\right) =\left\{ d_{3,0}^{i}d_{3,1}^{j}d_{3,2}^{k}\;|\ 0\leq i,j,k\leq
p-1\right\} \cup $ \newline
$\left\{ d_{3,0}^{p}d_{3,1}^{i}d_{3,2}^{j}\;|0\leq i,j\leq p-1\right\} \cup
\left\{ d_{3,1}^{p}d_{3,2}^{i}\;|\ 0\leq i\leq p-1\right\} \cup \left\{
d_{3,2}^{p}\right\} $\newline
is a free module basis for 
%TCIMACRO{\TeXButton{TeX field}{$\mathbb F_p$}}%
%BeginExpansion
$\mathbb F_p$%
%EndExpansion
$(3,1)$ over $D_{4}$.
\end{lemma}

%TCIMACRO{\TeXButton{Proof}{\proof}}%
%BeginExpansion
\proof%
%EndExpansion
Because of remark \ref{FreeModule}, our statement follows directly from the
following relations and induction on the total degree of $%
d_{3,2}^{m_{2}}d_{3,1}^{m_{1}}d_{3,0}^{m_{0}}$:\newline
i) $d_{4,0}=d_{4,3}d_{3,0}-d_{3,0}d_{3,2}^{p}$ $\Rightarrow $ $%
d_{3,0}d_{3,2}^{p}=-d_{4,0}+d_{4,3}d_{3,0}$;\newline
ii) $d_{4,1}=d_{4,3}d_{3,1}+d_{3,0}^{p}-d_{3,1}d_{3,2}^{p}$ $\Rightarrow $\ $%
d_{3,1}d_{3,2}^{p}=-$ $d_{4,1}+d_{4,3}d_{3,1}+d_{3,0}^{p}$;\newline
iii) $d_{4,2}=d_{4,3}d_{3,2}+d_{3,1}^{p}-d_{3,2}^{p+1}$ $\Rightarrow
d_{3,2}^{p+1}=-d_{4,2}+d_{4,3}d_{3,2}+d_{3,1}^{p}$;\newline
iv) $-d_{3,0}^{p+1}=d_{4,0}d_{3,1}-d_{4,1}d_{3,0}$;\newline
v) $-d_{3,0}d_{3,1}^{p}=d_{4,0}d_{3,2}-d_{4,2}d_{3,0}$;\newline
vi) $-d_{3,1}^{p+1}=d_{4,1}d_{3,2}-d_{4,2}d_{3,1}-d_{3,0}^{p}d_{3,2}$.%
%TCIMACRO{\TeXButton{End Proof}{\endproof}}%
%BeginExpansion
\endproof%
%EndExpansion

For each $t$, $1\leq t\leq n-1$, we define the set of all $\left( n-t\right) 
$-tuples\ \ 
\begin{equation*}
\mathcal{M}(n-2,t)=\{M=(p,m_{t},...,m_{n-2})\ |\ 0\leq m_{i}\leq p-1\}
\end{equation*}%
and, for each $M\in \mathcal{M}(n-2,t)$\ we define 
\begin{equation*}
d_{n-1}^{M}=d_{n-1,t-1}^{p}d_{n-1,t}^{m_{t}}...d_{n-1,n-2}^{m_{n-2}}
\end{equation*}

\begin{theorem}
\label{mimic}We have 
\begin{equation*}
B_{D_{n}}\left( 
%TCIMACRO{\TeXButton{TeX field}{\mathbb F_p}}%
%BeginExpansion
\mathbb F_p%
%EndExpansion
(n-1,1)\right) =\tbigcup\limits_{t=1}^{n-1}\{d_{n-1}^{M}\,|\,M\in \mathcal{M}%
(n-2,t)\}
\end{equation*}%
\textit{as a free module basis for 
%TCIMACRO{\TeXButton{TeX field}{$\mathbb F_p$}}%
%BeginExpansion
$\mathbb F_p$%
%EndExpansion
}$(n-1,1)$\textit{\ over }$D_{n}$\textit{.\ \ }
\end{theorem}

%TCIMACRO{\TeXButton{Proof}{\proof} }%
%BeginExpansion
\proof
%EndExpansion
Because of remark \ref{FreeModule}, we only have to prove that the given set
is a generating set. We use induction on the total degree $|m|:=\tsum m_{i}$
of a typical monomial $d^{m}=\prod\limits_{i=0}^{n-2}d_{n-1,i}^{m_{i}}$. 
\newline
Let us recall our relations:\newline
$d_{n,i}=d_{n,n-1}d_{n-1,i}-d_{n-1,i}d_{n-1,n-2}^{p}+d_{n-1,i-1}^{p}$ $%
\Rightarrow $\ 
\begin{equation}
d_{n-1,i}d_{n-1,n-2}^{p}=-d_{n,i}+d_{n,n-1}d_{n-1,i}+d_{n-1,i-1}^{p}\ \text{%
for }0\leq i\leq n-2  \label{rel-1}
\end{equation}%
$%
d_{n,i}d_{n-1,j}-d_{n,j}d_{n-1,i}=d_{n-1,i-1}^{p}d_{n-1,j}-d_{n-1,i}d_{n-1,j-1}^{p} 
$ $\Rightarrow $\ 
\begin{equation}
d_{n-1,i}d_{n-1,j-1}^{p}=-d_{n,i}d_{n-1,j}+d_{n,j}d_{n-1,i}+d_{n-1,i-1}^{p}d_{n-1,j}\ 
\text{for }0\leq i\leq j-1.  \label{rel2}
\end{equation}
\begin{equation}
d_{n-1,0}^{p+1}=d_{n,1}d_{n-1,0}-d_{n,0}d_{n-1,1}  \label{rel-3}
\end{equation}%
Please note that relations (\ref{rel-1}) and (\ref{rel-3}) reduce the total
degree. On the other hand, relation (\ref{rel2}) does not, but it moves the
same type of degree to the left with respect to index $i$.

It is obvious, because of the types of the relations above, that no other
relation can be deduced from the ones given. Namely, any combination of
these ends up to the one given.

Let $d_{i}$ denote $d_{n-1,i}$ for simplicity.

Let $d^{m}=\prod\limits_{i=1}^{l}d_{s_{i}}^{m_{s_{i}}}$, $0\leq
s_{1}<...<s_{l}\leq n-2$, and $0<m_{s_{i}}$.\ Let $f=d^{m}/d_{s_{1}}$. Then $%
f=\sum d\left( i\right) f\left( i\right) $\ where $d\left( i\right) \in
D_{n} $ and $f\left( i\right) $\ is a basis element by induction. Let $%
g=f(i)d_{s_{1}}$\ and $f\left( i\right) =\prod d_{s_{i}^{\prime
}}^{m_{s_{i}^{\prime }}}$. Here $m_{s_{1}^{\prime }}\leq p$\ and $%
m_{s_{i}^{\prime }}<p$.

If $s_{1}<s_{1}^{\prime }$, then $g$\ is a basis element. If $%
s_{1}=s_{1}^{\prime }$\ and $m_{s_{1}^{\prime }}<p$, then $g$\ is a basis
element.\ 

Let $s_{1}=s_{t}^{\prime }$, $m_{s_{t}^{\prime }}=p$\ and $t$ maximal. Thus $%
g=g^{\prime }d_{s_{t}^{\prime }}^{p+1}$. \newline
i) If $s_{t}^{\prime }=0$ or $n-2$, then the total degree of the
decomposition according to relations (\ref{rel-1})\ and (\ref{rel-3}) is
strictly less than that of $g$. \newline
ii) Let $0<s_{t}^{\prime }<n-2$. According to relation (\ref{rel2}), $%
d_{s_{t}^{\prime }}^{p+1}=d_{s_{t}^{\prime }-1}^{p}d_{s_{t}^{\prime
}+1}+others$.

(*) Now we consider $g^{\prime }d_{s_{t}^{\prime }-1}^{p}d_{s_{t}^{\prime
}+1}$. Again by relation (\ref{rel2}), this element is either a basis
element or decomposes to $g^{\prime }d_{s_{t}^{\prime
}-2}^{p}d_{s_{t}^{\prime }}d_{s_{t}^{\prime }+1}+others$. After a finite
number of steps either a basis element is obtained plus others or $%
d_{0}^{k_{0}}d_{1}^{k_{1}}...d_{s_{t}^{\prime }+1}^{k_{s_{t}^{\prime }+1}}$\
such that $k_{0}>p$\ and \ $k_{i}<p$. Now relation (\ref{rel-3}) is in
order.\ 

Let $s_{1}=s_{t}^{\prime }$, $m_{s_{t}^{\prime }}=p-1$\ and $t>1$. In this
case $g$ has the form $g=...d_{s_{t-1}^{\prime }}d_{s_{t}^{\prime }}^{p}...$
and we proceed as in (*) above. \ 
%TCIMACRO{\TeXButton{End Proof}{\endproof}}%
%BeginExpansion
\endproof%
%EndExpansion

\begin{corollary}
\label{mainCorollar}i) $\func{Im}\left( res^{\ast }:H^{\ast }\left( \Sigma
_{p}\int \Sigma _{p^{n-1}}\right) \rightarrow H^{\ast }\left( V\right)
\right) $ is isomorphic to a free module over $D_{n}$ on \ 
\begin{equation*}
\left\{ 
\begin{array}{c}
\hat{M}_{1,0}\hat{L}_{1}^{\left( p-2\right) }\hat{h}_{1}^{(p-1)m},\hat{M}%
_{n,s_{1},...,s_{k}}\hat{L}_{n}^{\left( p-2\right) }d_{n,0}^{\left( \left[ 
\frac{k+1}{2}\right] -1\right) }\hat{h}_{1}^{(p-1)m}\ |\  \\ 
0\leq m<A_{1},k\leq n,1\leq s_{k},0\leq s_{1}<...<s_{k}\leq n-1%
\end{array}%
\right\}
\end{equation*}%
Here $A_{1}=p^{n-1}+...+p$.

ii) $\func{Im}\left( res^{\ast }:H^{\ast }\left( \Sigma _{p^{n-1}}\int
\Sigma _{p}\right) \rightarrow H^{\ast }\left( V\right) \right) $ is
isomorphic to a free module over $D_{n}$ on \ 
\begin{equation*}
\left\{ 
\begin{array}{c}
\hat{M}_{n,n-1}\hat{L}_{n}^{\left( p-2\right) }f,\hat{M}%
_{n-1,s_{1},...,s_{k}}\hat{L}_{n-1}^{\left( p-2\right) }d_{n-1,0}^{\left( %
\left[ \frac{k+1}{2}\right] -1\right) }g\ |\  \\ 
f,g\in B_{D_{n}}\left( 
%TCIMACRO{\TeXButton{TeX field}{\mathbb F_p}}%
%BeginExpansion
\mathbb F_p%
%EndExpansion
(n-1,1)\right) \text{, }k\leq n-1,0\leq s_{1}<...<s_{k}\leq n-1%
\end{array}%
\right\}
\end{equation*}
\end{corollary}

%TCIMACRO{\TeXButton{Proof}{\proof}}%
%BeginExpansion
\proof%
%EndExpansion
This is an application of proposition \ref{1,n-1}, Theorem \ref{mimic} and %
\ref{invar-Bock-Parab}, corollary \ref{Corol-Restr-Bock-Upper} and lemma \ref%
{relatio-M}.%
%TCIMACRO{\TeXButton{End Proof}{\endproof}}%
%BeginExpansion
\endproof%
%EndExpansion

\section{The transfer}

In the opposite direction of the restriction map, a map is defined called
the transfer for $H$ a subgroup of finite index in $G$:%
\begin{equation*}
tr^{\ast }:H^{\ast }\left( H\right) \rightarrow H^{\ast }\left( G\right)
\end{equation*}%
At the cochain level $tr^{\ast }\left( a\right) \left( \lambda \right)
=\tsum\limits_{1}^{[G:H]}g_{i}a\left( g_{i}^{-1}\lambda \right) $. Here $%
a\in C_{H}^{i}=Hom_{%
%TCIMACRO{\TeXButton{TeX field}{\mathbb Z}}%
%BeginExpansion
\mathbb Z%
%EndExpansion
(H)}\left( C_{i},%
%TCIMACRO{\TeXButton{TeX field}{\mathbb F_p}}%
%BeginExpansion
\mathbb F_p%
%EndExpansion
\right) $, $\lambda \in C_{i}$\ and $\left\{ g_{i}\right\} $\ is a set of
left coset representatives (\cite{Ad-Mil} page 71).\ \ 

Let us recall from the introduction that the Weyl subgroups of $V$ in $%
\Sigma _{p^{n_{1}}}\int \Sigma _{p^{n_{2}}}$ and $\Sigma _{p^{n}}$ are $%
P(n_{1},n_{2})$ and the general linear group $GL(n,\mathbb{F}_{p})$
respectively. The induced inclusion 
\begin{equation*}
W_{\Sigma _{p^{n_{1}}}\int \Sigma _{p^{n_{2}}}}(V)\rightarrow W_{\Sigma
_{p^{n}}}\left( V\right)
\end{equation*}
induces 
\begin{equation*}
H^{\ast }\left( V\right) ^{P(n_{1},n_{2})}\overset{\tau ^{\ast }}{%
\rightarrow }H^{\ast }\left( V\right) ^{GL(n,\mathbb{F}_{p})}
\end{equation*}
given by $\tau ^{\ast }\left( f\right) =\tsum\limits_{1}^{[G:H]}g_{i}f$.
Here $f$ is a $P(n_{1},n_{2})$-invariant polynomial. $V$ is a 
%TCIMACRO{\TeXButton{TeX field}{$\mathbb F_p$}}%
%BeginExpansion
$\mathbb F_p$%
%EndExpansion
$G$-module. In our case the transfer is surjective and $H^{\ast }\left(
V\right) ^{GL(n,\mathbb{F}_{p})}$ is a direct summand. The following diagram
is commutative, please see \cite{Kuhn}.

\begin{equation*}
\begin{array}{ccc}
H^{\ast }\left( \Sigma _{p^{n_{1}}}\int \Sigma _{p^{n_{2}}}\right) & \overset%
{tr^{\ast }}{\rightarrow } & H^{\ast }\left( \Sigma _{p^{n}}\right) \\ 
\downarrow \left( res_{V}^{\Sigma _{p^{n_{1}}}\int \Sigma
_{p^{n_{2}}}}\right) ^{\ast } &  & \downarrow \left( res_{V}^{\Sigma
_{p^{n}}}\right) ^{\ast } \\ 
H^{\ast }\left( V\right) ^{W_{\Sigma _{p^{n_{1}}}\int \Sigma
_{p^{n_{2}}}}\left( V\right) } & \overset{\tau ^{\ast }}{\rightarrow } & 
H^{\ast }\left( V\right) ^{W_{\Sigma _{p^{n}}}\left( V\right) }%
\end{array}%
\end{equation*}

Campbell and Hughes (\cite{C-H}) have studied the transfer for the case: 
\begin{equation*}
\tau ^{\ast }:H_{n}\rightarrow D_{n}
\end{equation*}%
We extended the result above (\cite{Kech1}) for 
\begin{equation*}
\begin{array}{ccc}
H^{\ast }\left( V\right) ^{U_{n}} & \overset{\tau ^{\ast }}{\rightarrow } & 
H^{\ast }\left( V\right) ^{GL(n,\mathbb{F}_{p})}%
\end{array}%
\end{equation*}%
In this work we consider the induced map 
\begin{equation*}
\bar{\tau}^{\ast }:\func{Im}\left( res^{\ast }:H^{\ast }\left( G\right)
\rightarrow H^{\ast }\left( V\right) \right) \rightarrow \func{Im}\left(
res^{\ast }:H^{\ast }\left( \Sigma _{p^{n}}\right) \rightarrow H^{\ast
}\left( V\right) \right)
\end{equation*}%
Here $G=\Sigma _{p^{n},p}$, $\Sigma _{p}\int \Sigma _{p^{n-1}}$ and $\Sigma
_{p^{n-1}}\int \Sigma _{p}$.

Next we define a set of coset representatives for the groups under
consideration. We apply the method of Campbell and Hughes.

Let $\Pr_{n}\left( x\right) \in 
%TCIMACRO{\TeXButton{TeX field}{\mathbb F_p}}%
%BeginExpansion
\mathbb F_p%
%EndExpansion
[x]$ be an irreducible polynomial of degree $n\ $and$\ \sigma _{n}\ $ a root
of $\Pr_{n}\left( x\right) $\ in the $(p^{n}-1)$-st cyclotomic field $%
%TCIMACRO{\TeXButton{TeX field}{\mathbb F_p^n}}%
%BeginExpansion
\mathbb F_p^n%
%EndExpansion
$ over $%
%TCIMACRO{\TeXButton{TeX field}{\mathbb F_p}}%
%BeginExpansion
\mathbb F_p%
%EndExpansion
$ (\cite{Lidl}).Let $\sigma _{n}$ be a primitive root of unity and its
minimal polynomial 
\begin{equation*}
\Pr_{n}\left( x\right) =c_{0}+c_{1}x+...+c_{n-1}x^{n-1}+x^{n}
\end{equation*}%
Here there exists $j$ such that $c_{j}c_{0}\neq 0$. \ \ Then 
\begin{equation*}
\sigma _{n}^{n}=-\left( c_{0}+c_{1}\sigma _{n}+...+c_{n-1}\sigma
_{n}^{n-1}\right)
\end{equation*}%
and the companion matrix of $\Pr_{n}\left( x\right) $ is 
\begin{equation*}
A_{n}=\left( 
\begin{array}{cccc}
0 & \ldots & 0 & -c_{0} \\ 
1 &  &  & -c_{1} \\ 
\vdots & \ddots &  & \vdots \\ 
0 & \ldots & 1 & -c_{n-1}%
\end{array}%
\right)
\end{equation*}%
with $\Pr_{n}\left( A_{n}\right) =0_{n\times n}$. So $A_{n}$\ is a
representative for $\left( 
%TCIMACRO{\TeXButton{TeX field}{\mathbb F_p^n}}%
%BeginExpansion
\mathbb F_p^n%
%EndExpansion
\right) ^{\ast }$\ and can be identified with $\sigma _{n}$. 
\begin{eqnarray*}
%TCIMACRO{\TeXButton{TeX field}{\mathbb F_p^n} }%
%BeginExpansion
\mathbb F_p^n
%EndExpansion
&=&<\sigma _{n}^{0},\sigma _{n},...,\sigma _{n}^{n-1}>=<\sigma
_{n},...,\sigma _{n}^{p^{n-1}}>,<\sigma _{n}>=\left( 
%TCIMACRO{\TeXButton{TeX field}{\mathbb F_p^n}}%
%BeginExpansion
\mathbb F_p^n%
%EndExpansion
\right) ^{\ast } \\
&<&\sigma _{n}^{\frac{p^{n}-1}{p-1}}>=%
%TCIMACRO{\TeXButton{TeX field}{\mathbb F_p}}%
%BeginExpansion
\mathbb F_p%
%EndExpansion
^{\ast }
\end{eqnarray*}

$A_{n}$ acts linearly on $%
%TCIMACRO{\TeXButton{TeX field}{\mathbb F_p^n}}%
%BeginExpansion
\mathbb F_p^n%
%EndExpansion
$\ and $A_{n}\left( c\sigma _{n}^{i}\right) =c\sigma _{n}^{i+1}$. Let us
note that this action is compatible with the given action on the rings of
invariants: let $c\sigma _{n}^{i}$\ be represented by $\left(
0,...,0,c,0,...,0\right) $ with respect to the given basis, then $%
A_{n}\left( c\sigma _{n}^{i}\right) $ is the matrix multiplication between $%
\left( 0,...,0,c,0,...,0\right) $ and the $i+1$-th column of $A_{n}$. \ \ 

Moreover, $\sigma _{n}^{k}=\sigma _{n}^{k-n+1}\sigma _{n}^{n-1}$ or $%
A_{n}^{k-n+1}\left( \sigma _{n}^{n-1}\right) $\ for any $k$. Thus the last
column of $A_{n}^{k}$\ can be any non-zero element of $\left( 
%TCIMACRO{\TeXButton{TeX field}{\mathbb F_p^n}}%
%BeginExpansion
\mathbb F_p^n%
%EndExpansion
\right) ^{\ast }$.

Let 
\begin{equation*}
\Phi _{n}:<\sigma _{n}^{0},\sigma _{n},...,\sigma _{n}^{n-1}>\rightarrow
V^{n}
\end{equation*}%
Then the map induced by $\Phi _{n}(\sigma _{n}^{i})=y_{n-i}$ and linearity
is an isomorphism. Moreover, $A_{n}$ acts on $V\ $via $\Phi _{n}$: $%
A_{n}y_{i}=\Phi _{n}\left( \sigma _{n}\sigma _{n}^{n-i}\right) $. Now, $%
\left( 
%TCIMACRO{\TeXButton{TeX field}{\mathbb F_p^n}}%
%BeginExpansion
\mathbb F_p^n%
%EndExpansion
\right) ^{\ast }$ can be viewed as a subset of the group of automorphisms $%
GL(n,\mathbb{F}_{p})$.

Inductively we define $\sigma _{m}$ such that $\left\langle \sigma
_{m}\right\rangle \cong \left\langle y_{n-m+1},...,y_{n}\right\rangle ^{\ast
}$ and $\Phi _{m}:<\sigma _{m}^{0},\sigma _{m},...,\sigma _{m}^{m-1}>$ $%
\rightarrow \left\langle y_{n-m+1},...,y_{n}\right\rangle $. We consider $%
A_{m}\in GL(n,\mathbb{F}_{p})$ such that $A_{m}\left( y_{j}\right) =y_{j}$
for $1\leq j\leq n-m$\ and $A_{m}\left( y_{j}\right) =\Phi _{m}\left( \sigma
_{m}\sigma _{m}^{n-m+1-j}\right) $ for $1+n-m\leq j\leq n$. \ 

\begin{lemma}
i) Let \c{C}$=\{\left( A_{n}^{i}\right) ^{-1}\ |\ 0\leq i\leq p^{n}-2\}/\sim 
$\ where $A_{n}^{i}\sim A_{n}^{j}$, if there exists $c\in \left( 
%TCIMACRO{\TeXButton{TeX field}{\mathbb F_p}}%
%BeginExpansion
\mathbb F_p%
%EndExpansion
\right) ^{\ast }$\ such that $A_{n}^{i}=cA_{n}^{j}$.\ \ \ Then the set \c{C}
is a set of left coset representatives for $GL(n,\mathbb{F}_{p})$ over $%
P(1,n-1)$.

ii) Let \'{C}$=\{\left( A_{n}^{i}\right) ^{t}\ |\ 0\leq i\leq p^{n}-2\}/\sim 
$\ where $A_{n}^{i}\sim A_{n}^{j}$, if there exists $c\in 
%TCIMACRO{\TeXButton{TeX field}{\mathbb F_p}}%
%BeginExpansion
\mathbb F_p%
%EndExpansion
^{\ast }$\ such that $A_{n}^{i}=cA_{n}^{j}$.\ \ \ Then the set \'{C} is a
set of left coset representatives for $GL(n,\mathbb{F}_{p})$ over $P(n-1,1)$.
\end{lemma}

%TCIMACRO{\TeXButton{Proof}{\proof}}%
%BeginExpansion
\proof%
%EndExpansion
We recall that $|GL_{n}:P(1,n-1)|=\frac{p^{n}-1}{p-1}=|GL_{n}:P(n-1,1)|=|%
\text{\c{C}}|=|\acute{C}|$. The first column of $\left( A_{n}^{k}\right)
^{-1}$\ or the last row of $\left( A_{n}^{k}\right) ^{t}$ can be any
non-zero element of $\left( 
%TCIMACRO{\TeXButton{TeX field}{\mathbb F_p^n}}%
%BeginExpansion
\mathbb F_p^n%
%EndExpansion
\right) ^{\ast }$. Let $g,g^{\prime }\in P(1,n-1)$ and $\left(
A_{n}^{k}\right) ^{-1}g=\left( A_{n}^{l}\right) ^{-1}g^{\prime }$, then $%
\left( A_{n}^{k-l}\right) ^{-1}\in P(1,n-1)$\ which is not the case for $%
k\neq l$. The same is true for $P(n-1,1)$.%
%TCIMACRO{\TeXButton{End Proof}{\endproof}}%
%BeginExpansion
\endproof%
%EndExpansion
$\allowbreak $

The following proposition has been proved by Campbell and Hughes in \cite%
{C-H}.

\begin{proposition}
\label{coset,repr}The set $\left\{ \left( A_{n}^{i_{n}}\right)
^{-1}...\left( A_{1}^{i_{1}}\right) ^{-1}\ |\ 0\leq i_{m}\leq
p^{m}-2\right\} $ is a set of left coset representatives for $GL(n,\mathbb{F}%
_{p})$ over $U_{n}$.
\end{proposition}

%TCIMACRO{\TeXButton{Proof}{\proof}}%
%BeginExpansion
\proof%
%EndExpansion
We apply induction on $n$.%
%TCIMACRO{\TeXButton{End Proof}{\endproof}}%
%BeginExpansion
\endproof%
%EndExpansion

\begin{proposition}
Let $\xi :$ 
%TCIMACRO{\TeXButton{TeX field}{$\mathbb F_p$}}%
%BeginExpansion
$\mathbb F_p$%
%EndExpansion
$(1,n-1)\longrightarrow D_{n}$ be the natural epimorphism with respect to
the given free module basis $B$ and 
\begin{equation*}
\tau ^{\ast }:%
%TCIMACRO{\TeXButton{TeX field}{\mathbb F_p}}%
%BeginExpansion
\mathbb F_p%
%EndExpansion
(1,n-1)\rightarrow D_{n}
\end{equation*}
the transfer map. Then $\xi =\tau ^{\ast }$.
\end{proposition}

%TCIMACRO{\TeXButton{Proof}{\proof}}%
%BeginExpansion
\proof%
%EndExpansion
Let us recall that a free module basis consists of $\left(
h_{1}^{p-1}\right) ^{m}$\ for $0\leq m\leq \frac{p^{n}-p}{p-1}=p^{n-1}+...+p$%
. 
\begin{eqnarray*}
\tau ^{\ast }\left( h_{1}^{p-1}\right) ^{m} &=&\tsum_{i}\left(
A_{n}^{i}\right) ^{-1}\left( h_{1}^{p-1}\right) ^{m}=\tsum_{i}\left( \left(
A_{n}^{i}\right) ^{-1}h_{1}\right) ^{(p-1)m}= \\
\tsum_{u\in V}\left( u\right) ^{(p-1)m} &=&\left( p-1\right) \tsum_{1\leq
i\leq n,v\in <y_{1},...,y_{i-1}>}\left( y_{i}+v\right) ^{(p-1)m}
\end{eqnarray*}%
The last summand is a $GL$-invariant and so only $m\left( p-1\right)
=p^{n}-p^{k}$ for $1\leq k\leq n-1$\ should be considered, i.e. $\tau ^{\ast
}\left( h_{1}^{p-1}\right) ^{m}$ is a scalar multiple of $d_{n,k}$. Because
of proposition \ref{formula}, $d_{n,k}$ contains $\left(
\tprod\limits_{t=1}^{n-k}y_{k+t}^{p^{k+t-1}}\right) ^{p-1}$. Next we
consider the coefficient of this monomial in $\tau ^{\ast }\left(
h_{1}^{p-1}\right) ^{m}$\ or in\ \ 
\begin{equation*}
\tsum\limits_{v\in <y_{k+1},...,y_{n-1}>,u\in <y_{1},...,y_{k}>}\left(
y_{n}+v+u\right) ^{(p-1)m}
\end{equation*}%
This coefficient is \ $p^{k}\frac{\left( p^{n}-p^{k}\right) !}{%
\tprod\limits_{t=1}^{n-k}\left( p^{k+t-1}\left( p-1\right) \right) !}\equiv 0%
\func{mod}p$. Thus $\tau ^{\ast }\left( h_{1}^{p-1}\right) ^{m}=0$.\ \ 
%TCIMACRO{\TeXButton{End Proof}{\endproof}}%
%BeginExpansion
\endproof%
%EndExpansion

\begin{remark}
According to the last proof, if $m\left( p-1\right) =p^{n}-1$, then 
\begin{equation*}
\tau ^{\ast }\left( h_{1}^{p-1}\right) ^{m}=\left( p-1\right) d_{n,0}
\end{equation*}
\end{remark}

\begin{theorem}
\cite{C-H} Let $\xi :$ $H_{n}\longrightarrow D_{n}$ be the natural
epimorphism with respect to the given basis $B$ and $\tau ^{\ast
}:H_{n}\rightarrow D_{n}$ the transfer map. Then $\xi =\tau ^{\ast }$.
\end{theorem}

%TCIMACRO{\TeXButton{Proof}{\proof}}%
%BeginExpansion
\proof%
%EndExpansion
It is obvious that $\tprod \left( A_{s}^{i_{s}}\right) ^{-1}\left( \tprod
h_{t}^{r_{t}}\right) =$ 
\begin{equation*}
\left( A_{n}^{i_{n}}\right) ^{-1}(h_{1}^{r_{1}}\left(
A_{n-1}^{i_{n-1}}\right) ^{-1}(h_{2}^{r_{2}}...\left( A_{2}^{i_{2}}\right)
^{-1}(h_{n-1}^{r_{n-1}}(\left( A_{1}^{i_{1}}\right) ^{-1}h_{n}^{r_{n}}))...)
\end{equation*}%
Let $r_{i}<p^{n-i+1}-1$, then the proof of last proposition for $n=n-i+1$
implies that $\tsum\limits_{i_{n-i+1}}\left( A_{n-i+1}^{i_{n-i+1}}\right)
^{-1}h_{i}^{r_{i}}=0$. Now the statement follows. 
%TCIMACRO{\TeXButton{End Proof}{\endproof}}%
%BeginExpansion
\endproof%
%EndExpansion

\begin{corollary}
Let $\xi :$ $%
%TCIMACRO{\TeXButton{TeX field}{\mathbb F_p}}%
%BeginExpansion
\mathbb F_p%
%EndExpansion
(1,...,1)\longrightarrow D_{n}$ be the natural epimorphism with respect to
the given basis $B$ and $\tau ^{\ast }:%
%TCIMACRO{\TeXButton{TeX field}{\mathbb F_p}}%
%BeginExpansion
\mathbb F_p%
%EndExpansion
(1,...,1)\rightarrow D_{n}$ the transfer map. Then $\xi =\tau ^{\ast }$.
\end{corollary}

Next we consider $P\left( n-1,1\right) $. In this case the use of the coset
representatives arises technical problems. Instead, using degree arguments,
we shall prove that only particular elements of the given basis might be
expressed with respect to Dickson algebra generators. Then applying Steenrod
operations on Dickson algebra generators, we shall prove that the transfer
map $\tau ^{\ast }:%
%TCIMACRO{\TeXButton{TeX field}{\mathbb F_p}}%
%BeginExpansion
\mathbb F_p%
%EndExpansion
(n-1,1)\rightarrow D_{n}$ coincides with the natural epimorphism $\xi :$ 
%TCIMACRO{\TeXButton{TeX field}{$\mathbb F_p$}}%
%BeginExpansion
$\mathbb F_p$%
%EndExpansion
$(n-1,1)\longrightarrow D_{n}$ with respect to basis $B$.

The next technical lemma will be needed for the proof of our next Theorem.

\begin{lemma}
Let $m_{i}$ and $m_{j}^{\prime }$ be non-negative integers such that $0\leq
m_{i}\leq p-1$, $m_{j}^{\prime }\geq 0$ and $p$ a prime number.

1) Let $0\leq i\leq n-2$ and $0\leq j\leq n-1$. Then the equation 
\begin{equation*}
\tsum\limits_{0}^{n-2}m_{i}\left( p^{n-1}-p^{i}\right)
=\tsum\limits_{0}^{n-1}m_{j}^{\prime }\left( p^{n}-p^{j}\right)
\end{equation*}%
does not have an integral solution.

2) Let $i_{0}<i\leq n-2$ and $0\leq j\leq n-1$. Then the equation 
\begin{equation*}
\left( p^{n}-p^{i_{0}+1}\right) +\tsum\limits_{i_{0}+1}^{n-2}m_{i}\left(
p^{n-1}-p^{i}\right) =\tsum\limits_{0}^{n-1}m_{j}^{\prime }\left(
p^{n}-p^{j}\right)
\end{equation*}%
admits solutions of type $m_{i_{0}+1}=...=m_{k}=p-1$, $m_{i}=0$ for $k<i$
and $m_{n-1}^{\prime }=k-i_{0}$, $m_{k+1}^{\prime }=1$ and zero otherwise. \
\ \ Here $i_{0}<k\leq n-2$ and $m_{i}=0$ for any $i>i_{0}$, $%
m_{i_{0}+1}^{\prime }=1$ and zero otherwise.
\end{lemma}

\begin{theorem}
\label{Main}Let $\xi :$ 
%TCIMACRO{\TeXButton{TeX field}{$\mathbb F_p$}}%
%BeginExpansion
$\mathbb F_p$%
%EndExpansion
$(n-1,1)\longrightarrow D_{n}$ be the natural epimorphism with respect to
the given free module basis $B$ and $\tau ^{\ast }:%
%TCIMACRO{\TeXButton{TeX field}{\mathbb F_p}}%
%BeginExpansion
\mathbb F_p%
%EndExpansion
(n-1,1)\rightarrow D_{n}$ the transfer map. Then $\xi =\tau ^{\ast }$.
\end{theorem}

%TCIMACRO{\TeXButton{Proof}{\proof} }%
%BeginExpansion
\proof
%EndExpansion
If we show that $\tau ^{\ast }(d)=0$ for all $d$ in the basis, then $\xi
=\tau ^{\ast }$. Because of the statement in last lemma only the following
cases should be considered: $d_{n-1,i_{0}}^{p}$, and $%
d_{n-1,i_{0}}^{p}d_{n-1,i_{0}+1}^{p-1}...d_{n-1,k}^{p-1}$.

Let $\tau ^{\ast }\left( d_{n-1,i_{0}}^{p}\right) =cd_{n,i_{0}+1}$. Applying 
$P^{p^{i_{0}}}$, we get $\tau ^{\ast }\left( d_{n-1,i_{0}-1}^{p}\right)
=cd_{n,i_{0}}$. Applying $P^{p}...P^{p^{i_{0}}-1}$ on the previous element,
we get 
\begin{equation*}
\tau ^{\ast }\left( d_{n-1,0}^{p}\right) =cd_{n,1}
\end{equation*}%
But 
\begin{equation*}
P^{1}\tau ^{\ast }\left( d_{n-1,0}^{p}\right) =0\neq P^{1}\left(
cd_{n,1}\right) =cd_{n,0}
\end{equation*}

Let $\tau ^{\ast }\left(
d_{n-1,i_{0}}^{p}d_{n-1,i_{0}+1}^{p-1}...d_{n-1,k}^{p-1}\right)
=cd_{n,k+1}d_{n,n-1}^{k-i_{0}}$. \newline
Let $f=d_{n-1,i_{0}}^{p}d_{n-1,i_{0}+1}^{p-1}...d_{n-1,k}^{p-1}$ and $%
g=d_{n,k+1}d_{n,n-1}^{k-i_{0}}$.\ We would like to apply $P^{p^{k}}$. Using
Theorem \ref{Steenr-action-dn} we show that no monomial of $P^{p^{k}}f\neq 0$
is in the ideal $\left( d_{n-1,k}\right) $. Then $P^{p^{k}}f\in \left(
d_{n-1,0},...,d_{n-1,k-1}\right) $. Since 
\begin{equation*}
p^{k}=\left( p-1\right) \left( p^{k-1}+...+p^{i_{0}}\right) +p^{i_{0}}
\end{equation*}
$P^{p^{k}}f$\ contains the summand $%
d_{n-1,i_{0}-1}^{p}d_{n-1,i_{0}}^{p-1}...d_{n-1,k-1}^{p-1}$. So $%
P^{p^{k}}f\neq 0$.\ 

Claim: No monomial of $P^{p^{k}}f$ is in the ideal $\left( d_{n-1,k}\right) $%
. We prove the claim by showing that there does not exist a solution of \ \ 
\begin{equation*}
p^{k}=\tsum\limits_{j=i_{0}}^{k-1}\tsum\limits_{i=0}^{j}a_{j,i}p^{i}
\end{equation*}%
unless $a_{k-1,k-1}=p-1$. In that case 
\begin{equation*}
P^{p^{k}}f=\tsum\limits_{\left( m_{0},...,m_{k-1}\right) }c\left(
m_{0},...,m_{k-1}\right) \tprod\limits_{t=0}^{k-1}d_{n-1,t}^{m_{t}}
\end{equation*}%
Here $0\leq a_{j,i}\leq a_{j,i+1}\leq p-1$ for $i_{0}+1\leq j\leq k-1$\ and $%
0\leq a_{i_{0},i}\leq a_{i_{0},i+1}\leq 1$. We consider the extreme cases
and prove that there is no positive solution.

Let $a_{k-1,k-1}=p-2$, $a_{k-1,i}=p-2$, $a_{j,i}=p-1$ for\ $i_{0}+1\leq
j\leq k-2$\ and $a_{i_{0},t}=1$. Then \ \ \ \ 
\begin{equation*}
\tsum\limits_{j=i_{0}}^{k-1}\tsum\limits_{i=0}^{j}a_{j,i}p^{i}=p^{k}-\left(
k-i_{0}-3\right) <p^{k}
\end{equation*}%
Now the claim follows.

It is obvious that if $a_{k-1,k-1}=p-1$, then no summand of $P^{p^{k}}f$\ is
in $\left( d_{n-1,k}\right) $.

Applying $P^{1}...P^{p^{k}}$, we get%
\begin{equation*}
P^{1}...P^{p^{k}}f=0\text{ and }P^{1}...P^{p^{k}}\left(
d_{n,k}d_{n,n-1}^{k-i_{0}}\right) =d_{n,0}d_{n,n-1}^{k-i_{0}}
\end{equation*}%
The last line is an application of Theorem \ref{Steenr-action-dn}. 
%TCIMACRO{\TeXButton{End Proof}{\endproof}}%
%BeginExpansion
\endproof%
%EndExpansion

The next example is a counterexample to the statement of last Theorem in the
case $\func{Im}\left( res_{V}^{\Sigma _{p^{n},p}}\right) ^{\ast }\supsetneq
H_{n}$.

\begin{example}
\label{Paradeigma}Let $p=3$\ and $n=2$. Then $\Pr_{2}\left( x\right)
=2+x+x^{2}$ \ and $\sigma _{2}^{2}=2\sigma _{2}+1$ or $A_{2}=\left( 
\begin{array}{cc}
0 & 1 \\ 
1 & 2%
\end{array}%
\right) $, $A_{2}^{-1}=\left( 
\begin{array}{cc}
1 & 1 \\ 
1 & 0%
\end{array}%
\right) $. A set of coset representatives for $GL(2,3)$\ over $P(1,1)$\ is
given in proposition \ref{coset,repr}.

By direct computation, $\tau ^{\ast }\left( M_{1,0}h_{1}^{p^{2}-1-p}\right)
=M_{2,1}L_{2}^{p-2}\neq 0$. Let us note that $M_{1,0}h_{1}^{p^{2}-1-p}$\ is
a basis element.\ 
\end{example}

\begin{theorem}
\cite{Kech1}\label{final}Let $\xi ,\bar{\tau}^{\ast }:\func{Im}\left(
res_{V}^{\Sigma _{p^{n},p}}\right) ^{\ast }\rightarrow \func{Im}\left(
res_{V}^{\Sigma _{p^{n}}}\right) ^{\ast }$. Then $\xi =\bar{\tau}^{\ast }$
in the ideal generated by $\left( d_{n,0}\right) $. \ 
\end{theorem}

%TCIMACRO{\TeXButton{Proof}{\proof} }%
%BeginExpansion
\proof
%EndExpansion
Let $f\in \func{Im}\left( res_{V}^{\Sigma _{p^{n},p}}\right) ^{\ast }$. $\ \ 
$Then $\bar{\tau}^{\ast }\left( fd_{n,0}\right) =\bar{\tau}^{\ast }\left(
f\right) d_{n,0}$. But according to lemma \ref{relatio-M}: 
\begin{equation*}
\bar{\tau}^{\ast }\left( fd_{n,0}\right) =\bar{\tau}^{\ast }\left(
\tsum\limits_{J}M_{n;J}L_{n}^{p-2}h^{I\left( J\right) }\right)
=\tsum\limits_{J}M_{n;J}L_{n}^{p-2}\bar{\tau}^{\ast }\left( h^{I\left(
J\right) }\right) =
\end{equation*}%
\begin{equation*}
\tsum\limits_{J}M_{n;J}L_{n}^{p-2}\xi \left( h^{I\left( J\right) }\right)
\end{equation*}%
If $\xi \left( h^{I\left( J\right) }\right) $ is not divisible by $d_{n,0}$,
then it must be zero. \ 
%TCIMACRO{\TeXButton{End Proof}{\endproof}}%
%BeginExpansion
\endproof%
%EndExpansion

\end{document}